\documentclass[mnsc,nonblindrev]{informs3_hide}
\RequirePackage{tgtermes}
\RequirePackage{bm}
\RequirePackage{endnotes}

\OneAndAHalfSpacedXII %

\usepackage{algorithm}
\usepackage{algpseudocode}
\usepackage{tikz}

\usepackage{dsfont}
\usepackage{pifont}
\usepackage{amsmath,amsfonts,amssymb}
\usepackage{mathtools}
    \allowdisplaybreaks
\usepackage{bm}
\usepackage[mathscr]{euscript}
\usepackage[inline]{enumitem}
\usepackage{xcolor}
\usepackage{graphicx}
\usepackage{multirow}
\usepackage{bigstrut}
\usepackage{booktabs}
\usepackage{microtype}
\usepackage{fix-cm}
\usepackage{physics}
\usepackage{xfrac}
\usepackage{nicematrix}
\usepackage{caption}
\usepackage[hidelinks]{hyperref}

\DeclareMathOperator{\ind}{\mathds 1}
\DeclareMathOperator{\pr}{\mathbb P}
\DeclareMathOperator{\E}{\mathbb E}
\DeclareMathOperator{\N}{\mathbb N}

\DeclareMathOperator{\vari}{\mathbb V\mathrm{ar}}

\DeclareMathOperator{\ci}{\mathsf{CI}}
\DeclareMathOperator{\diag}{\mathrm{diag}}

\newcommand{\GGG}{\mathbb G}
\newcommand{\QQ}{\mathbb Q}
\newcommand{\RR}{\mathbb R}

\newcommand{\msfV}{\mathsf V}

\newcommand{\scrA}{\mathcal A}

\newcommand{\scrD}{\mathcal D}

\newcommand{\scrF}{\mathcal F}

\newcommand{\scrI}{\mathcal I}

\newcommand{\scrL}{\mathcal L}

\newcommand{\scrS}{\mathcal S}

\newcommand{\scrV}{\mathcal V}

\newcommand{\scrX}{\mathcal X}
\newcommand{\scrY}{\mathcal Y}
\newcommand{\scrZ}{\mathcal Z}

\usepackage{subfigure}

\usepackage{soul}
\usepackage{comment}
\definecolor{myRed}{RGB}{139,0,0}

\usepackage{natbib}
 \bibpunct[, ]{(}{)}{,}{a}{}{,}%
 \def\bibfont{\small}%

\EquationsNumberedThrough    %

\TheoremsNumberedThrough     %
\ECRepeatTheorems  %

\MANUSCRIPTNO{}

\begin{document}

\RUNAUTHOR{Wang and Zhang}

\RUNTITLE{``Over-optimizing'' for Normality}

\TITLE{``Over-optimizing'' for Normality: Budget-constrained Uncertainty Quantification for Contextual Decision-making}

\ARTICLEAUTHORS{%
\AUTHOR{Yanyuan Wang, Xiaowei Zhang}
\AFF{Department of Industrial Engineering and Decision Analytics, The Hong Kong University of Science and Technology, Clear Water Bay, Hong Kong SAR, \EMAIL{yanyuan.wang@connect.ust.hk}, \EMAIL{xiaoweiz@ust.hk}}
} %

\ABSTRACT{We study uncertainty quantification for contextual stochastic optimization, focusing on weighted sample average approximation (wSAA), which uses machine-learned relevance weights based on covariates.
Although wSAA is widely used for contextual decisions, its uncertainty quantification remains limited. In addition, computational budgets tie sample size to optimization accuracy, creating a coupling that standard analyses often ignore.
We establish central limit theorems for wSAA and construct asymptotic-normality-based confidence intervals for optimal conditional expected costs. We analyze the statistical--computational tradeoff under a computational budget, characterizing how to allocate resources between sample size and optimization iterations to balance statistical and optimization errors.
These allocation rules depend on structural parameters of the objective; misspecifying them can break the asymptotic optimality of the wSAA estimator.
We show that ``over-optimizing'' (running more iterations than the nominal rule) mitigates this misspecification and preserves asymptotic normality, at the expense of a slight slowdown in the convergence rate of the budget-constrained estimator.
The common intuition that ``more data is better'' can fail under computational constraints: increasing the sample size may worsen statistical inference by forcing fewer algorithm iterations and larger optimization error.
Our framework provides a principled way to quantify uncertainty for contextual decision-making under computational constraints.
It offers practical guidance on allocating limited resources between data acquisition and optimization effort, clarifying when to prioritize additional optimization iterations over more data to ensure valid confidence intervals for conditional performance.
}


\maketitle

\section{Introduction} \label{sec:intro}

Contextual stochastic optimization (CSO) has recently attracted substantial attention in operations research and management science, driven by applications in uncertain, data-rich decision-making environments such as inventory management \citep{BertsimasKallus20}, portfolio optimization \citep{ElmachtoubGrigas22}, and service capacity management \citep{NotzPibernik22}. In CSO, environmental uncertainty is represented by a random outcome (e.g., product demand) whose distribution is unknown and depends on observable covariates (e.g., search volume as a proxy for consumer attention). Using historical observations of covariates and outcomes, the decision-maker seeks to minimize expected cost under the \emph{conditional} distribution of the random outcome given a new covariate observation.
The resulting decision is tailored to the current context rather than applied uniformly across all contexts, allowing the use of predictive information to improve operational performance.

A wide range of methods have been proposed for CSO, integrating machine learning tools---linear and kernel regression, tree-based methods, and neural networks---to model either the conditional expected cost or a policy that maps covariates to decisions \citep{BanRudin19,BertsimasKallus20,QiGrigasShen21,BertsimasKoduri22,HoNguyenKilincKarzan22,KannanBayraksanLuedtke25}.
These methods produce a data-driven decision for a new context, effectively a \emph{point estimate} of the context-specific optimal decision.

Performance analyses of these methods typically establish asymptotic optimality of this point estimate: as the historical data grows, its expected cost converges to that of the true optimal decision. However, point estimates alone can be risky because they provide no measure of the uncertainty from random variation in historical data.
The issue is exacerbated in contextual settings: conditioning on covariates to customize decisions increases variability relative to the non-contextual case, and the uncertainty can grow as more covariates are used for refined customization.
Accurate \emph{uncertainty quantification} is therefore essential to complement CSO methods and to support reliable, risk-aware decision-making.

Despite this need, uncertainty quantification has been underexplored in the CSO literature.
Equally crucial for practice---but often overlooked---are \emph{computational constraints}.
CSO methods lead to data-driven optimization problems. Most analyses study statistical properties of the optimizer while implicitly assuming an oracle solver.
In reality, these problems are solved by iterative algorithms whose computational expense scales with the size of the historical data.
Under a fixed time budget, using a larger sample may reduce the number of iterations the algorithm can run, and the algorithm may stop before reaching the true optimum.
As a result, the quality of the solution obtainable within a given time frame depends on sample size, directly linking statistical choices (how much data to use) to computational feasibility and, ultimately, to the reliability of uncertainty quantification.

To illustrate, consider the widely used CSO method of weighted sample average approximation (wSAA) \citep{BertsimasKallus20}, which is the focus of our paper.
This method estimates the conditional expected cost of a decision when the outcome distribution varies with covariates. It proceeds in two steps: (i) evaluate costs on historical outcomes, and (ii) form a weighted sum of these costs, where each weight reflects how relevant the associated covariate is to the new covariate. This weighting approximates the conditional distribution of the outcome given the new covariate. Solving the resulting \emph{wSAA problem}, which minimizes the weighted estimator of the conditional expected cost, yields an estimate of the optimal decision for the new context.

When solving the wSAA problem with iterative methods (e.g., gradient descent), the per-iteration computational expense scales with the size of the historical data because computing the gradient requires traversing all observations.
Under a fixed computational budget, the sample size therefore limits the number of iterations, introducing a statistical--computational tradeoff that is pivotal to budget allocation.

\begin{figure}[ht]
    \FIGURE{\includegraphics[width=\textwidth]{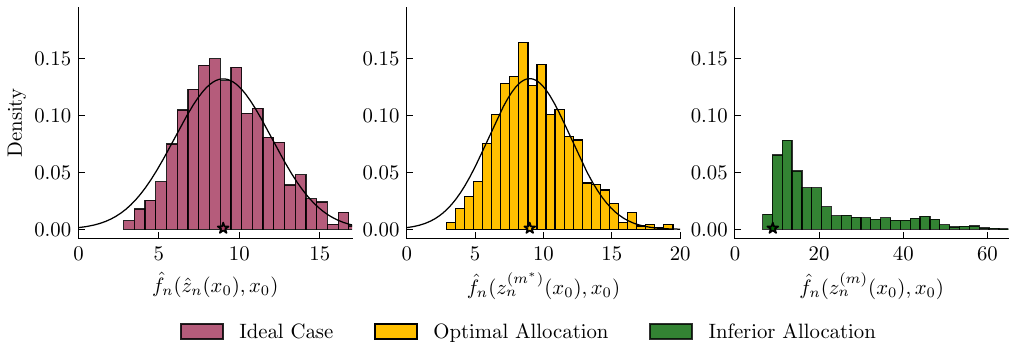}}
{Distribution of wSAA Estimates of Optimal Conditional Expected Cost. \label{fig:moti}}
{Stars mark the true optimal conditional expected cost for a newsvendor problem (see Section~\ref{app:exp-nv} of the e-companion).
The histograms are based on 1,000 replications, and the overlaid density curves show the asymptotic normal distribution implied by our theory.
The left panel presents an idealized setting in which the wSAA problem (sample size $10^4$) is solved to optimality (Section~\ref{subsec:tradeoff-ideal}).
The middle panel shows a budget-constrained setting with a total budget of $10^5$, optimally allocated between sample size and iterations (Section~\ref{subsec:tradeoff-constr}).
The right panel illustrates how a na{\"i}ve treatment of the statistical--computational tradeoff leads to inferior budget allocation and flawed asymptotic-normality-based uncertainty quantification.}
\end{figure}

In this setting, the usual intuition that ``more data is better'' can fail. Using all historical data---or choosing the sample size without accounting for per-iteration expense---can increase finite-sample variability and undermine uncertainty quantification. Figure~\ref{fig:moti} shows how na{\"i}ve budget allocation can induce a complex finite-sample distribution for the wSAA estimator of the optimal conditional expected cost,
which invalidates asymptotic-normality-based confidence intervals and causes uncertainty quantification to hinder rather than support contextual decision-making.

\subsection{Main Contributions}

Our first contribution is to establish central limit theorems (CLTs) for wSAA when the weights are constructed nonparametrically via Nadaraya–Watson kernel regression.
These CLTs permit the construction of confidence intervals for the conditional expected costs of optimal decisions in CSO problems.
In contrast to SAA \citep{Shapiro21}, the terms in the wSAA objective are dependent, even when historical covariates are independent.
The dependence arises because each weight is built from the full dataset, inducing correlations across terms. As a result, standard SAA analysis does not carry over directly to wSAA.
We address the interdependent weights by combining kernel regression techniques with empirical process theory \citep{Pollard90}. We show that the wSAA estimator of the conditional expected cost---as a function of the decision---converges, after proper scaling, to a Gaussian process. This functional limit result yields CLTs for optimal conditional expected costs.

Second, we characterize a statistical--computational tradeoff in uncertainty quantification under computational constraints. Using more historical data reduces the statistical error in wSAA but raises the computational expense of each algorithm iteration.
With a fixed budget, higher per-iteration expense means fewer iterations, which increases optimization error and yields a \emph{biased} estimate of the optimal conditional expected cost.
To address this tradeoff, we derive CLTs for a \emph{budget-constrained wSAA estimator}, where the bias diminishes asymptotically as the computational budget grows.
The convergence rates depend on the optimization algorithm’s convergence regime (linear, sublinear, or superlinear) so that the optimization error is small relative to the statistical error, which is necessary to eliminate asymptotic bias.

Third, we show the benefits of over-optimizing the wSAA problem to make uncertainty quantification more robust.
Optimal budget allocation relies on structural parameters of the objective (e.g., Lipschitz constant).
These problem-specific parameters are typically unknown, and misspecifying them can lead to allocation choices that break the asymptotic normality of the budget-constrained wSAA estimator, which is crucial for constructing confidence intervals based on normal limits.
Running a few more iterations than the budget-optimal choice offsets such misspecification and preserves asymptotic normality, at the cost of a slight reduction in the convergence rate of the budget-constrained estimator.
In this sense, over-optimizing---prioritizing normality over the fastest convergence---yields more reliable uncertainty quantification and is reassuring for practitioners.

\subsection{Related Works}

The CSO literature has expanded substantially in recent years, with diverse methodological developments and applications; see \cite{QiShen22} and \cite{Sadana25} for comprehensive surveys. As a leading approach in this field, the wSAA method provides a versatile framework for solving CSO problems by integrating various machine learning techniques to compute weights that measure the relevance of historical covariate observations to new contextual information.
Early work by \cite{Hannah10} and \cite{BanRudin19} explored the use of Nadaraya-Watson kernel regression for weight construction, while \cite{BertsimasKallus20} broadened the methodology to include $k$-nearest neighbors, local linear regression, decision trees, and random forests. \cite{KallusMao23} further enhanced tree-based weight assignment  by integrating optimization objectives  into the training of tree-based models, moving beyond conventional supervised learning approaches that only rely on historical data.

The wSAA method has been applied beyond traditional CSO problems. \cite{NotzPibernik22} and \cite{BertsimasMcCordSturt23}  extended the framework to multi-stage settings, and \cite{RahimianPagnoncelli23} considered chance-constrained problems. High-dimensional covariates can degrade performance by encouraging overfitting in the models used for weight construction. To address this issue, several approaches have emerged. For example, \cite{Srivastava21} and \cite{LinChenLiShen22}  introduced an additional regularization term into the objective function.

Existing analysis of the wSAA method has primarily focused on establishing asymptotic optimality or deriving generalization error bounds
\citep{BanRudin19,BertsimasKallus20,Srivastava21}.
While those bounds indicate the method's average performance over the distribution of the covariates, they do not quantify uncertainty in the estimates of optimal \emph{conditional} expected costs specific to new observations.
In contrast, our CLTs enable the construction of confidence intervals for optimal conditional expected costs with asymptotically exact coverage guarantees.

Our theoretical analysis contributes to uncertainty quantification for contextual decision-making, an area that remains underexplored except for several recent studies. \cite{CaoGaoYang21} examined statistical inference for parametric models in CSO, which does not cover wSAA due to its nonparametric nature.
The parametric models either describe the relationship between covariates and the random variable involved in the cost function, or directly model the relationship between covariates and the conditional expected cost.
Instead of establishing CLTs as we do, they derived concentration bounds on the parameter estimates, which in turn yield confidence regions for the true parameter values.
\cite{Cao24} leveraged conformal prediction techniques \citep{AngelopoulosBates23} to create confidence intervals for optimal decisions of a contextual newsvendor problem in which the policy class may be misspecified. However, this approach relies heavily on the unique connection between quantile regression and the newsvendor cost function and does not readily generalize to other cost functions.
\cite{Iyengar24} investigated the use of cross-validation techniques for uncertainty quantification in supervised learning and CSO problems.
However, their focus is to construct confidence intervals for \emph{unconditional} expected costs of optimal policies---marginalized over the distribution of covariates,
whereas our focus is on optimal \emph{conditional} expected costs.
In addition, none of these studies analyze the statistical--computational tradeoff involved in uncertainty quantification under computational budgets.

\subsection{Notation and Organization }

The following notation is used throughout the paper.
All vectors $v$ are treated as column vectors, with $v^\intercal$ denoting the transpose and $\|v\| \coloneqq \sqrt{v^\intercal v}$ denoting the Euclidean norm. For a matrix $A$, $\|A\|$ denotes its spectral norm.
For a positive integer $n$,
$[n]$ denotes the index set $\{ 1, \ldots, n\}$.
For a real-valued sequence $a_n$, we say
$a_n \asymp b_n$ if there exist constants $C',C > 0$ such that $C' < a_n/b_n < C$ when $n$ is large enough, and $a_n \sim b_n$ if $a_n/b_n \to 1$ as $n \to \infty$.
We say $X_n = o_{\pr}(1)$ if for any $\epsilon>0$, $\pr(|X_n|>\epsilon) \to 0$ as $n \to \infty$, and $X_n = O_{\pr}(1)$ if $X_n$ is \emph{stochastically bounded}; that is, for any $\epsilon > 0$, there exists a constant $M > 0$ such that $\pr(|X_n| > M) < \epsilon$ for $n$ large enough.
Suppose $Y_n$ is another sequence of random variables. We say $X_n = o_{\pr}(Y_n)$ and $X_n = O_{\pr}(Y_n)$ if $X_n/Y_n = o_{\pr}(1)$ and $X_n/Y_n = O_{\pr}(1)$, respectively.

The remainder of this paper is organized as follows.
Section~\ref{sec:form} formulates the CSO problem and presents its wSAA counterpart.
Section~\ref{sec:bdgt-alloc} studies uncertainty quantification for wSAA, first in an idealized setting without computational constraints and then under computational budgets, where the wSAA problem is solved by a linearly convergent optimization algorithm.
Section~\ref{sec:overopt} discusses the benefits of over-optimizing to preserve asymptotic normality for reliable uncertainty quantification, especially when convergence parameters may be misspecified.
Section~\ref{sec:extension} extends the statistical--computational tradeoff analysis to sublinear and superlinear optimization algorithms.
Section~\ref{sec:ci} provides practically implementable confidence intervals.
Section~\ref{sec:exper} validates our theoretical results through numerical experiments.
Section~\ref{sec:concl} concludes with remarks on future studies.
Omitted proofs and additional numerical results appear in the e-companion.

\section{Problem Formulation}\label{sec:form}

A decision-maker seeks a cost-minimizing decision where the cost function $F(z;Y)$ depends on both the decision variable $z $ and a random variable $Y$.
While the value of $Y$ is not yet realized when the decision must be made, its distribution depends on observable covariates $X$ that represent the available contextual information.
The decision-making problem can thus be expressed as
\begin{align} \label{eq:cso}
\min_{z \in \scrZ} \Bigg\{ f(z,x_0) \coloneqq \E[F(z;Y)|X = x_0] \Bigg\}, %
\end{align}
where $\scrZ \subseteq \RR^{d_z}$ represents the feasible set, the expectation is taken with respect to the conditional distribution of $Y$ given $X=x_0$, and $f(z,x_0)$  denotes the conditional expected cost.

In practice, the conditional distribution  is usually unknown, but the decision-maker may have access to a dataset $\scrD_n = \{(x_i,y_i)\}_{i=1}^n$ containing $n$ independent and identically distributed (i.i.d.) historical observations of $(X,Y)$.
For simplicity, we assume $X $ and $Y $ follow continuous distributions on their respective supports $\scrX \subseteq \RR^{d_x}$ and $\scrY \subseteq \RR^{d_y}$, although our theory extends to cases with discrete or mixed-valued $X$.
We also assume that the uncertainty is exogenous (i.e., the decision variable $z$ does not affect the distribution of $Y$).
The wSAA method can adapt to endogenous uncertainty when historical decisions are available and fully determined by the covariates (i.e., no unobserved confounders, or $Y$ and $Z$ are independent conditioning on $X$); see \cite{BertsimasKallus20}.

\subsection{Weighted Sample Average Approximation} \label{subsec:wsaa}

The wSAA method for solving the CSO problem~\eqref{eq:cso} approximates the objective function $f(z,x_0)$ with a weighted sample average of the cost function evaluated on the dataset $\scrD_n$:
\begin{align} \label{eq:cso-wSAA}
\min_{z \in \scrZ} \Bigg\{ \hat{f}_n(z,x_0) \coloneqq \sum_{i=1}^n w_n(x_i, x_0) F(z;y_i) \Bigg\}, %
\end{align}
where the weights $w_n(x_i, x_0)$ satisfy $\sum_{i=1}^n w_n(x_i, x_0) = 1$ and $w_n(x_i, x_0) \ge 0$ for all $i \in [n]$, measuring the relevance of each historical observation $x_i$ to the new observation $x_0$.
In essence, the wSAA method approximates the conditional distribution of $Y$ given $X=x_0$ by an empirical distribution where each outcome $y_i$ appears with probability $w_n(x_i,x_0)$.
When all weights equal $1/n$, the approximated conditional expected cost function $\hat{f}_n(z,x_0)$ reduces to the sample average of $\{F(z, y_i)\}_{i=1}^n$, and the wSAA problem~\eqref{eq:cso-wSAA} becomes an SAA problem that ignores the dependence of $Y$ on $X$ and simply makes decisions without contextual information.

In this paper, we focus on weights constructed using Nadaraya-Watson kernel regression, a leading approach in the literature \citep{BanRudin19,BertsimasKallus20}.
Specifically,
\begin{align} \label{eq:nw-weight}
w_{n}(x_i,x_0) = \frac{K((x_i-x_0) /h_n)}{\sum_{j=1}^n K((x_j-x_0) /h_n)},%
\end{align}
where $K:\RR^{d_x}\mapsto \RR$ is a kernel function, and $h_n>0$ is the bandwidth parameter which usually takes the form of $h_n=h_0 n^{-\delta}$ for some $h_0 > 0$ and $\delta \in (0,1)$.
Common kernel functions include the uniform kernel $K(u) = \ind(\|u\|\le 1)$,
Epanechnikov kernel $K(u) = (1-\|u\|^2) \ind(\|u\|\le 1)$,
and Gaussian kernel $K(u) = \exp(-\|u\|^2/2)$,
where $\ind(\cdot)$ denotes the indicator function.

Let $z^*(x_0)$ and $\hat{z}_n(x_0)$ denote optimal solutions to the CSO problem~\eqref{eq:cso} and the wSAA problem~\eqref{eq:cso-wSAA}, respectively.
Under mild conditions on the cost function $F$, kernel function $K$, and bandwidth $h_n$, \cite{BertsimasKallus20} showed that $\hat{z}_n(x_0)$ is asymptotically optimal, with the optimality gap in conditional expected cost vanishing; that is, $f(\hat{z}_n(x_0), x_0) - f(z^*(x_0), x_0) \to  0$ as $n\to\infty$ almost surely.
Specifically for newsvendor problems,
\cite{BanRudin19} established generalization bounds on the gap between the true value $f(z^*(x_0), x_0)$ and its estimate $\hat{f}_n(\hat{z}_n(x_0), x_0)$, averaged over the distribution of $x_0$.

Our main goal in this paper is to quantify uncertainty in the estimates of $f(z^*(x_0), x_0)$ obtained from solving the wSAA problem.
We consider two scenarios,
depending on whether computational expenses of solving this problem are taken into account.
In Section~\ref{subsec:tradeoff-ideal}, assuming the wSAA problem can be solved to optimality, which yields the solution $\hat{z}_n(x_0)$, we derive a CLT for the estimate $\hat{f}_n(\hat{z}_n(x_0), x_0)$ for any given $x_0$.
Based on this result, we can construct confidence intervals for $f(z^*(x_0), x_0)$ with asymptotically exact coverage guarantees.
In Section \ref{subsec:tradeoff-constr} and beyond, we extend our theory to factor in computational expenses of solving the wSAA problem via an iterative optimization algorithm.
With a finite budget, the algorithm returns an approximated solution $z^{(m)}_n(x_0)$ after $m$ iterations.
We establish CLTs for the \emph{budget-constrained} estimate $\hat{f}_n(z^{(m)}_n(x_0), x_0)$, facilitating uncertainty quantification under computational budgets.
Since each iteration's computational expense is proportional to $n$, the convergence rates of these CLTs reflect a tradeoff between statistical error (decreasing with $n$) and optimization error (decreasing with $m$). We analyze three types of optimization algorithms---linearly, sublinearly, and superlinearly convergent---which affect such a tradeoff through their respective convergence behaviors.

\begin{remark}
In addition to kernel regression, other machine learning methods can be used to construct the weight function $w_n(\cdot, x_0)$.
Our theory extends to weights constructed using $k$-nearest neighbors or local linear regression, due to their close relationship with kernel regression \citep{PaganUllah99}.
We omit this extension as it would introduce substantial technical complexity without offering significant new insights.
Developing CLTs for weights constructed using tree-based models \citep{BertsimasKallus20,KallusMao23}---potentially under computational budgets---requires a different approach
and remains beyond our scope.
\end{remark}

\subsection{Basic Assumptions}\label{subsec:assumptions}

\begin{assumption}[Regularity] \label{assump:regularity}
    \begin{enumerate}[label=(\roman*)]
        \item $F(z;y)$ is Lipschitz continuous in $z \in \scrZ$ uniformly for $y \in \scrY$:  there exists a constant $C_F>0$ such that $|F(z;y)-F(z';y)| \le C_F \|z-z'\|$ for all $y\in \scrY$ and $z, z'\in \scrZ$.
        \label{assump:F-lipschitz}
        \item $\scrZ $ is a nonempty and compact subset of $\RR^{d_z}$. \label{assump:Z-compact}

    \end{enumerate}
\end{assumption}

\begin{assumption}[Smoothness] \label{assump:smooth-1} \leavevmode
For any $x_0\in \scrX$, let $\scrV(x_0)$ be an open neighborhood of $x_0$.
    \begin{enumerate}[label=(\roman*)]
        \item For all $z \in \scrZ$, $f(z,x)$ is twice continuously differentiable in $x\in\scrX$. For any $x_0\in\scrX$, there exists a constant $C_f > 0$ such that  $\sup_{z \in \scrZ}|f(z,x)| \leq C_f$, $\sup_{z \in \scrZ}\|\nabla_x f(z,x)\| \leq C_f$ and $\sup_{z \in \scrZ}\|\nabla_x^2 f(z,x)\| \leq C_f$ for all $x\in \scrV(x_0)$. \label{assump:smooth-1-f}
        \item $X$ has a density $p(x)$ that is twice continuously differentiable and bounded away from zero.  For any $x_0\in\scrX$, there  exists a constant  $C_p > 0$  such that $|p(x)| \leq C_p$, $\|\nabla_x p(x)\| \leq C_p$ and $\|\nabla_x^2 p(x)\| \leq C_p$ for all $x\in \scrV(x_0)$. \label{assump:smooth-1-p}
        \item For all $z,z' \in \scrZ$, $\nu(z,z',x) \coloneqq \E[F(z;Y)F(z';Y)|X=x]$ is continuously differentiable in $x\in\scrX$.  For any $x_0\in\scrX$, there  exists a constant $C_{\nu} > 0$ such that $\sup_{z,z' \in \scrZ}|\nu(z,z',x)| \leq C_{\nu}$ and $\sup_{z,z' \in \scrZ}\|\nabla_x \nu(z,z',x)\| \leq C_{\nu}$ for all  $x\in \scrV(x_0)$. \label{assump:smooth-1-nu}
    \end{enumerate}
\end{assumption}

\begin{assumption}[Envelope] \label{assump:envelope}
    There exists an envelope function $M: \scrY \mapsto \RR_+$ for the class of functions $\scrF \coloneqq \{F(z;\cdot): z \in \scrZ\}$  such that
    \begin{enumerate}[label=(\roman*)]
        \item $\sup_{z \in \scrZ} |F(z;y)| \le M(y) < \infty$ for almost every $y \in \scrY$, and
        \item $\E[M^{2+\gamma}(Y)|X=x] \le C_M$ for all $x\in\scrV(x_0)$ and some $\gamma > 0$, where $C_M > 0$ is a constant.
    \end{enumerate}
\end{assumption}

\begin{assumption}[Kernel] \label{assump:kern}
    \begin{enumerate}[label=(\roman*)]
        \item The kernel function $K:\RR^{d_x}\mapsto \RR_+$ is spherically symmetric and has a finite second-order moment:  $\int_{\RR^{d_x}} u K(u) du = 0$ and $\int_{\RR^{d_x}} uu^\intercal K(u) du = \Upsilon(K) I_{d_x}  $ for some positive constant $\Upsilon(K) < \infty$, where $I_{d_x}$ denotes the $d_x \times d_x$  identity matrix.
        \item $\sup_{u \in \RR^{d_x}} |K(u)| < \infty$.
        \item $\|u\|^{d_x} K(u) \to 0$ as $\|u\| \to \infty$.
        \item
        $\int_{\RR^{d_x}} |K(u)|^{2+\gamma} du< \infty$ for some $\gamma > 0$.

 \label{assump:kern-integral}
    \end{enumerate}
\end{assumption}

\begin{assumption}[Bandwidth] \label{assumption:band}
    For all $n\geq 1$, $h_n=h_0 n^{-\delta}$, for some constants $h_0 > 0$ and $\delta \in (1/(d_x+4), 1/d_x)$.
\end{assumption}

\begin{assumption}[Uniqueness]\label{assump:uni-sol}
The CSO problem~\eqref{eq:cso} has a unique optimal solution $z^*(x_0)$.
\end{assumption}

Assumption~\ref{assump:regularity} rules out pathological cases.
Part~\ref{assump:F-lipschitz} can be relaxed to two conditions: \begin{enumerate*}[label=(\alph*)] \item $|F(z;y)-F(z';y)| \le L(y) \|z-z'\|$ for all $y\in \scrY$ and $z, z'\in \scrZ$ with $\E[L^{2+\gamma}(Y)|X=x_0] \le C_F$ for some constant $C_F>0$; and
\item $F(z;y)$ is equicontinuous in $z \in \scrZ$.
\end{enumerate*}
For part~\ref{assump:Z-compact}, the compactness of $\scrZ$ can be relaxed to
three conditions: \begin{enumerate*}[label=(\alph*)] \item $\scrZ$ is closed; \item $\liminf_{\|z\| \to \infty} \inf_{y \in \scrY} F(z;y) > -\infty$; and \item for any $x\in\scrX$, there exists $S_x\subset \scrY$ such that $\pr(y \in S_x \vert X = x)>0$ and $F(\cdot;y)$ is uniformly coercive over $S_x$; that is, $\lim_{\|z\|\to\infty} F(z;y) = \infty$ uniformly over $y \in S_x$. \end{enumerate*} These relaxed conditions parallel the Weierstrass theorem in deterministic optimization \citep[Theorem 3.2.1]{Bertsekas16}. They ensure $f(\cdot,x_0)$ and $\hat{f}_n(\cdot,x_0)$ are coercive. Combined with the equicontinuity of $F(\cdot;y)$ and boundedness of $\scrZ$, this guarantees finite, attainable optimal values $f(z^*(x_0), x_0)$ and $\hat{f}_n(\hat{z}_n(x_0), x_0)$.

Assumption~\ref{assump:smooth-1} ensures that $f(z,\cdot)$, $p(\cdot)$, and $\nu(z,z',\cdot)$ are sufficiently smooth near $x_0$, analogous to standard assumptions in the asymptotic analysis of nonparametric regression estimators.
While the first- and second-order differentiability conditions could be relaxed to Lipschitz and H\"older continuity, respectively, we maintain these stronger conditions for expositional clarity.

Assumption~\ref{assump:envelope} posits an envelope function to control the values of functions in the class  $\scrF$ and excludes infinite values of $F(z;\cdot)$ on null sets. The finite conditional moment condition on $M(Y)$ of order higher than two holds when the conditional distribution of $Y$ has sufficiently fast tail decay, as observed in sub-Gaussian and sub-exponential cases.
This is analogous to  Lyapunov's condition for classic CLTs \citep[Section 27]{Billingsley95}.

Assumption~\ref{assump:kern} stipulates basic kernel function properties.
Common kernels introduced in Section~\ref{subsec:wsaa} (e.g., Gaussian) readily satisfy this assumption.
Assumptions~\ref{assump:envelope} and~\ref{assump:kern} parallel the standard regularity conditions used in the asymptotic normality theory of kernel regression.

Assumption~\ref{assumption:band} follows the standard bandwidth choice, with one key difference: we require $\delta > 1/(d_x + 4)$ instead of the usual $\delta \in (0, 1/d_x)$.
This stronger condition serves to debias the estimator $\hat{f}_n(z, x_0)$
when establishing a functional central limit theorem (FCLT) for the objective function estimates, as elaborated in Section~\ref{sec:bdgt-alloc}.

Assumption~\ref{assump:uni-sol} states that the optimal solution set of the CSO problem is a singleton.
When multiple optimal solutions exist,
the expectation of the infimum of a limiting Gaussian process over such a set is typcially nonzero.
The nonlinearity introduced by the infimum operator complicates the construction of confidence intervals.
In contrast, when the optimal solution is unique, the infimum operator becomes unbinding: it reduces to evaluating the limiting process at this optimal point. As a result,
$\hat{f}_n(\hat{z}_n(x_0),x_0)$ is asymptotically normal, as established in Theorem \ref{thm:clt-cont}.

\section{Statistical--Computational Tradeoffs} \label{sec:bdgt-alloc}
We first consider a idealized scenario without computational constraints to investigate how the uncertainty underlying the wSAA estimates relates to the sample size (Section~\ref{subsec:tradeoff-ideal}). This prepares us to analyze  the statistical--computational tradeoff in the budget-constrained case (Section~\ref{subsec:tradeoff-constr}).

\subsection{Idealized Case: No Computational Constraints} \label{subsec:tradeoff-ideal}

To derive a CLT for the wSAA estimator of the optimal conditional expected cost, we proceed in two steps.
First, we establish an FCLT for $\hat{f}_n(\cdot,x_0)$.
Second, we apply the delta method with a first-order expansion of the min-value function---following \cite{Shapiro89,Shapiro91}'s analysis of directional derivatives---to obtain the desired CLT for $\hat{f}_n(\hat{z}_n(x_0), x_0)$.

One special consideration in deriving the FCLT for the wSAA estimator $\hat{f}_n(\cdot, x_0)$ concerns the dependence among weights $\{w_n(x_i, x_0): i\in[n]\}$, which arises because each of them is computed using the same dataset $\scrD_n$ as defined in \eqref{eq:nw-weight}.
To resolve this, we combine empirical process theory  with  kernel regression analysis.
This approach potentially extends to cases where observations in $\scrD_n$ are not i.i.d. but generated from a strong mixing process, following \cite{AndrewsPollard94}.

Specifically, we write $\hat{f}_n(z, x_0) = \hat{r}_n(z, x_0) / \hat{p}_n(x_0)$ where
\[
\hat{r}_n(z,x_0) \coloneqq \frac{1}{nh_n^{d_x}} \sum_{i=1}^n K\left( \frac{x_i-x_0}{h_n}\right) F(z;y_i)\quad\mbox{and}\quad
\hat{p}_n(x_0) \coloneqq \frac{1}{nh_n^{d_x}} \sum_{i=1}^n K\left( \frac{x_i-x_0}{h_n}\right),
\]
and write $f(z, x_0) = r(z, x_0) / p(x_0)$ where
$r(z,x_0) \coloneqq \int_{\scrY} F(z;y) p(x_0,y) d y$
with $p(x,y)$ denoting the joint density of $X$ and $Y$.
Using Taylor's expansion, we obtain
\begin{equation}\label{eq:taylor}
\begin{aligned}
        \hat{f}_n(z,x_0) - f(z,x_0)
        = {}& \frac{\hat{r}_n(z,x_0)  - f(z,x_0) \hat{p}_{n}(x_0)}{p(x_0)} \\
        & + O_{\pr}\left(| \hat{r}_n(z,x_0) -  r(z,x_0)| |\hat{p}_n(x_0) - p(x_0)| + |\hat{p}_n(x_0) - p(x_0)|^2\right).
\end{aligned}
\end{equation}
The first term in \eqref{eq:taylor} equals
\begin{align}\label{eq:taylor-term-1}
    \frac{1}{p(x_0)} \sum_{i=1}^n  K\left(\frac{x_i-x_0}{h_n}\right)(F(z;y_i)-f(z, x_0)),
\end{align}
representing a scaled difference between the expectation of $g(z; (X, Y))$ under the empirical distribution of $(X,Y)$ and that under the true distribution. Here,
$g(z; (x,y)) \coloneqq  K((x - x_0)/h_n) F(z; y) / p(x_0)$.
This formulation permits the use of empirical process theory
\citep{Pollard90} to establish an FCLT from~\eqref{eq:taylor-term-1}.
To analyze the higher-order terms in \eqref{eq:taylor},
we examine the errors $\hat{r}_n(z,x_0) - r(z,x_0)$ and
$\hat{p}_n(x_0) - p(x_0)$ separately.
For the former,
an FCLT can be derived using empirical process theory, similar to the analysis in \eqref{eq:taylor-term-1}.
For the latter,
the kernel density estimator $\hat{p}_n(x_0)$ has well-established asymptotic properties \citep[Theorem~2.10]{PaganUllah99_ec}. Combining the analyses of the leading term and higher-order terms in \eqref{eq:taylor} leads to Proposition~\ref{prop:fclt-cont}.

\begin{proposition} \label{prop:fclt-cont}
Under Assumptions~\ref{assump:regularity}--\ref{assumption:band},
    \[n^{(1-\delta d_x)/2} (\hat{f}_n(\cdot,x_0) - f(\cdot,x_0)) \Rightarrow h_0^{-d_x/2} \GGG(\cdot, x_0),\]
    as $n\to\infty$,
    for all $x_0 \in \scrX$, where $\GGG(\cdot, x_0)$ is a zero-mean Gaussian process with covariance function
    \[
    \Psi(z,z',x_0) = \frac{R_2(K)}{p(x_0)} \E[(F(z;Y) - f(z,x_0)) (F(z';Y) - f(z',x_0))|X=x_0].
    \]
    Here, $R_2(K) \coloneqq \int_{\RR^{d_x}} K^2(u) du < \infty$.
\end{proposition}

Note that this FCLT holds for all $\delta\in(1/(d_x+4), 1/d_x)$ (Assumption~\ref{assumption:band}).
This naturally raises the question of optimal bandwidth selection.
In kernel regression literature, it is typically addressed via mean squared error (MSE) minimization.
For given $z$ and $x_0$,
$\hat{f}_n(z, x_0)$ has a bias of order $O(n^{-2\delta} + n^{-(1-\delta d_x)})$ and a variance of order $O(n^{-(1-\delta d_x)})$, where $\delta\in(0, 1/d_x)$ \citep[pp.~101--103]{PaganUllah99}.
Minimizing the MSE by equating the orders of the squared bias and the variance gives $\delta = 1/(d_x + 4)$.
Although analogous FCLTs can be derived for any $\delta \in (0,1/(d_x+4)]$,
we assume a narrower range of $\delta$ for the purpose of debiasing.
When $\delta \leq 1/(d_x + 4)$, the limiting Gaussian process in the corresponding FCLT has a nonzero mean, indicating a non-vanishing bias.
This bias involves derivatives of unknown quantities such as $\nabla_x f(z,x_0)$, $\nabla_x^2 f(z,x_0)$, and $\nabla_x p(x_0)$, which are difficult to estimate, rendering the construction of confidence intervals for $f(z, x_0)$ challenging.
To circumvent this issue, we choose $\delta$ to be strictly larger than its optimal value under the MSE criterion---for example, $\delta = 1/(d_x + 4 - \epsilon)$ with a small $\epsilon>0$.
While this choice slightly compromises the convergence rate, it removes the bias from $\hat{f}_n(z, x_0)$, thereby simplifying statistical inference.
This approach---known as ``undersmoothing''---is commonly used in kernel regression literature \citep{Hall92}.

With the FCLT for $\hat{f}_n(\cdot,x_0)$ in place, we follow the approach of \cite{Shapiro89,Shapiro91} to derive the CLT for the wSAA estimator of the optimal conditional expected cost.
We treat $f(\cdot, x_0)$ and $\hat{f}_n(\cdot, x_0)$ as random elements of $C(\scrZ)$, the Banach space of continuous functions $\phi:\scrZ \mapsto \RR$ endowed with the sup-norm $|\phi| = \sup_{z \in \scrZ}|\phi(z)|$.
For the min-value function $\vartheta(\phi) \coloneqq \inf_{z \in \scrZ} \phi(z)$,
its Hadamard directional derivative at $\varphi \in C(\scrZ)$ in the direction $\varsigma \in C(\scrZ)$ is given by $\vartheta'_{\varphi}(\varsigma) = \inf_{z \in \scrS^*(\varphi)}\varsigma(z)$, where $\scrS^*(\varphi) = \argmin_{z \in \scrZ} \varphi(z)$.
Applying the delta method \citep[Theorem~9.74]{Shapiro21_ec} with $\varphi = f(\cdot,x_0)$ and $\varsigma = \hat{f}_n(\cdot,x_0) - f(\cdot,x_0)$, we obtain
\begin{align*}
    \hat{f}_n(\hat{z}_n(x_0),x_0) - f^*(x_0)
    = \inf_{z \in \scrZ^*(x_0)} (\hat{f}_n(z,x_0) - f(z,x_0)) + o_{\pr}\left(n^{-(1-\delta d_x)/2}\right),
\end{align*}
where $\scrZ^*(x_0)$ denotes the set of optimal solutions to the CSO problem~\eqref{eq:cso}.
To simplify notation, we abbreviate $f(z^*(x_0),x_0)$ as $f^*(x_0)$ hereinafter, whenever there is no risk of confusion.
Combining this first-order expansion with the FCLT for $\hat{f}_n(\cdot,x_0)$ leads to
\[n^{(1-\delta d_x)/2} (\hat{f}_n(\hat{z}_n(x_0),x_0) - f^*(x_0)) \Rightarrow \inf_{z \in \scrZ^*(x_0)} \GGG(z, x_0).\]

\begin{theorem} \label{thm:clt-cont}
    Under Assumptions~\ref{assump:regularity}--\ref{assump:uni-sol},
    \begin{align*}
        n^{(1-\delta d_x)/2} (\hat{f}_n(\hat{z}_n(x_0),x_0) - f^*(x_0)) \Rightarrow h_0^{-d_x/2} N(0,\msfV(z^*(x_0),x_0)),
    \end{align*}
    as $n\to\infty$,
    for all $x_0 \in \scrX$,
    where $N(0,v)$ is a normal distribution with mean $0$ and variance $v$, and
    \begin{equation}\label{eq:V-expression}
    \msfV(z,x_0) \coloneqq \frac{\sigma^2(z,x_0)}{p(x_0)} R_2(K).
    \end{equation}
    Here, $\sigma^{2}(z,x) \coloneqq \E[(F(z;Y) - f(z, x))^{2}|X=x]$ is the conditional variance of $F(z;Y)$ given $X=x$.
\end{theorem}

\subsection{Uncertainty Quantification with Computational Constraints} \label{subsec:tradeoff-constr}
In real-world applications, where computational budgets are limited, the question of how much historical data to use becomes crucial.
While increasing the sample size $n$ improves the accuracy of the wSAA estimator at the rate specified in Theorem~\ref{thm:clt-cont}, it becomes  computationally more expensive to solve the wSAA problem~\eqref{eq:cso-wSAA}. In particular, the computational expenses of evaluating the objective function or its gradient scale with $n$.

Under computational constraints, the algorithm may be forced to terminate with a suboptimal solution. Thus, using all available data may not be optimal when computational resources are limited---particularly if the optimization algorithm requires many iterations to converge due to the objective function lacking certain structural properties.
In the following, we derive a CLT to characterize the statistical--computational tradeoff that arises when solving the wSAA problem with a given budget.

Let $\scrA$ denote the optimization algorithm used to solve the wSAA problem. Starting from an initial point $z_n^{(0)}(x_0)\in\scrZ$, it iteratively generates a sequence $\{z_n^{(m)}(x_0)\}_{m \in \N_+}$. (The choice of the initial point may depend on $n$ and $x_0$, which in turn affects all subsequent iterates.) Let $\Gamma$ denote the computational budget available to decision-makers.
Since the objective function $\hat{f}_n(\cdot, x_0)$ of the wSAA problem~\eqref{eq:cso-wSAA} consists of $n$ terms,
the computational expense for one single iteration of $\scrA$ (e.g., gradient descent) is typically proportional to $n$.
Without loss of generality, we assume that the computational expense per iteration is normalized to 1. Otherwise, one could rescale $\Gamma$ by this expense without affecting asymptotic analysis.
Therefore, the number of iterations $m$ permitted by the computational budget satisfies $nm \approx \Gamma$.

Under these computational budget constraints, the decision-maker obtains a solution $z^{(m)}_n(x_0)$ from the algorithm instead of the optimal solution $\hat{z}_n(x_0)$ to the wSAA problem.
The difference between the budget-constrained wSAA estimator $\hat{f}_n(z^{(m)}_n(x_0), x_0)$  and the optimal conditional expected cost $f^*(x_0)$ can be decomposed into two parts:
\begin{align}\label{eq:error-decomp}
    \hat{f}_n(z^{(m)}_n(x_0), x_0) - f^*(x_0) =
    \underbrace{\hat{f}_n(z^{(m)}_n(x_0), x_0) - \hat{f}_n(\hat{z}_n(x_0), x_0)}_{\text{optimization error}} + \underbrace{\hat{f}_n(\hat{z}_n(x_0), x_0) -  f^*(x_0)}_{\text{statistical error}}.
\end{align}

Theorem~\ref{thm:clt-cont} establishes that the statistical error is $O_{\pr}(n^{-(1-\delta d_x)/2})$, which diminishes as the sample size $n\to\infty$. On the other hand,  the optimization error diminishes as the number of iterations $m\to\infty$, at a rate that depends on both the optimization algorithm and the properties of the wSAA objective function. Consequently, the computational budget constraints give rise to a tradeoff between the two types of errors. We examine how to allocate a budget $\Gamma$ between $n$ and $m$ to minimize the total error in~\eqref{eq:error-decomp}. Particularly, we consider \emph{asymptotically admissible} allocation rules, meaning that as $\Gamma\to\infty$, both $n$ and $m$ grow without bound and $nm/\Gamma \to 1$.

In this subsection, we focus on the case where $\scrA$ exhibits linear convergence (defined below), and defer the discussion on algorithms of other types (sublinear and superlinear) to Section~\ref{sec:extension}.

\begin{definition}[Linear Convergence] \label{defi:lin}
    An algorithm is said to converge linearly for solving the wSAA problem~\eqref{eq:cso-wSAA} if there exists a constant $\theta \in (0,1)$ such that
    \begin{align}
        \hat{f}_n(z^{(m)}_n(x_0),x_0) - \hat{f}_n(\hat{z}_n(x_0),x_0) \le \theta (\hat{f}_n(z^{(m-1)}_n(x_0),x_0) - \hat{f}_n(\hat{z}_n(x_0),x_0)), \label{eq:linear-conv}
    \end{align}
    for all $z_n^{(0)}(x_0) \in \scrZ$ and $n, m \in \N_+$.
\end{definition}

\begin{example}\label{example:linear-conv}
    If $\hat{f}_n(\cdot, x_0)$ is $\lambda$-strongly convex and differentiable with $L$-Lipschitz continuous derivatives, then gradient descent with Armijo backtracking
    achieves linear convergence. Specifically, the constant $\theta$ can be taken as $1 - 4a(1 - a)b\lambda/L$, where $a \in (0,0.5)$ and $b \in (0,1)$ are the line search parameters \citep[Theorem~1.3.17]{Polak97}.  If a projection operator is applied in each iteration, $\theta$ becomes $1 - ab\lambda/L$ \citep[Theorem~1.3.18]{Polak97}.
Note that $\hat{f}_n(\cdot, x_0)$ satisfies this condition when $F(\cdot;y)$ is differentiable and $\lambda(y)$-strongly convex for almost every $y \in \scrY$, where $L = C_F$ under Assumption~\ref{assump:regularity}-\ref{assump:F-lipschitz} and $\lambda = \sum_{i=1}^n w_n(x_i,x_0) \lambda(y_i)$.
\end{example}

By recursively applying inequality \eqref{eq:linear-conv}, we can bound the optimization error by:
\begin{align}
    \hat{f}_n(z^{(m)}_n(x_0),x_0) - \hat{f}_n(\hat{z}_n(x_0),x_0) \leq \theta^m (\hat{f}_n(z^{(0)}_n (x_0),x_0) - \hat{f}_n(\hat{z}_n(x_0),x_0)) = O_{\pr}(\theta^m). \label{eq:opti-error-linear}
\end{align}
For a given budget $\Gamma$, minimizing the total error in~\eqref{eq:error-decomp} requires balancing the optimization error $O_{\pr}(\theta^m)$ with the statistical error $O_{\pr}(n^{-(1 - \delta d_x)/2})$. This leads to the asymptotically optimal budget allocation rule: $m(\Gamma) \sim \kappa \log(\Gamma)$ and $n(\Gamma) \sim \Gamma/(\kappa\log(\Gamma))$ for some constant $\kappa>0$.
Theorem~\ref{thm:lin} formalizes this intuition by establishing the asymptotic normality of $\hat{f}_n(z^{(m)}_n(x_0),x_0)$, provided that $\kappa$ exceeds a certain threshold.

\begin{theorem} \label{thm:lin}
Suppose Assumptions~\ref{assump:regularity}--\ref{assump:uni-sol} hold, and $\scrA$ is linearly convergent with parameter $\theta\in(0,1)$.
Consider an asymptotically admissible budget allocation $\{n(\Gamma),m(\Gamma)\}_{\Gamma \in \N_+}$ that satisfies $m(\Gamma) \sim \kappa \log \Gamma$ as $\Gamma\to\infty$ for some constant $\kappa > 0$.
    \begin{enumerate}[label=(\roman*)]
        \item If $\kappa \ge (1-\delta d_x)/(2\log (1/\theta))$,
        then
        \begin{align}
            \left(\frac{\Gamma}{\log\Gamma}\right)^{(1-\delta d_x)/2} (\hat{f}_n(z^{(m)}_n(x_0),x_0) - f^*(x_0)) \Rightarrow \left(\frac{\kappa^{1-\delta d_x}}{h_0^{d_x}}\right)^{1/2} N(0,\msfV(z^*(x_0),x_0)), \label{eq:thm-lin-i}
        \end{align}
        as $\Gamma\to\infty$,  for all $x_0 \in \scrX$,
        where $\msfV(z^*(x_0),x_0)$ is defined in \eqref{eq:V-expression}.
        \item If $0 < \kappa < (1-\delta d_x)/(2\log (1/\theta))$, then
        \begin{align}
            \Gamma^{\kappa \log (1/\theta)} (\hat{f}_n(z^{(m)}_n(x_0),x_0) - f^*(x_0)) = O_{\pr}(1), \label{eq:thm-lin-ii}
        \end{align}
        as $\Gamma\to\infty$, for all $x_0 \in \scrX$.
    \end{enumerate}
\end{theorem}

Let $\kappa^* = (1-\delta d_x)/(2\log (1/\theta))$. When $\kappa \geq \kappa^*$, the CLT in \eqref{eq:thm-lin-i} implies that the error of the budget-constrained wSAA estimator diminishes at a rate of $(\Gamma/\log\Gamma)^{-(1 - \delta d_x)/2}$. While this rate is independent of $\kappa$, the asymptotic variance is proportional to $\kappa^{(1-\delta d_x)/2}$. Since $\delta d_x < 1$, it decreases with $\kappa$. Therefore, among all choices of $\kappa \geq \kappa^*$, setting $\kappa = \kappa^*$ minimizes the asymptotic variance of the error $\hat{f}_n(z^{(m)}_n(x_0),x_0) - f^*(x_0)$.
When $\kappa < \kappa^*$, the error behaves as $O_{\pr}(\Gamma^{-\kappa \log(1/\theta)})$. Since $\theta\in(0,1)$, this rate is faster with a larger $\kappa$ and approaches $O_{\pr}(\Gamma^{- (1 - \delta d_x)/2 })$ as $\kappa$ increases, eventually matching the rate of $(\Gamma/\log\Gamma)^{-(1 - \delta d_x)/2}$ up to a logarithmic factor when $\kappa=\kappa^*$.
Analyzing both scenarios confirms that $\kappa^*$ is indeed optimal given $\delta$ and $\theta$. With $\kappa^*$, the budget-constrained wSAA estimator converges at a rate of $(\Gamma/\log\Gamma)^{-(1 - \delta d_x)/2}$, nearly recovering the rate $n^{-(1 - \delta d_x)/2}$ established in Theorem~\ref{thm:clt-cont},
with only a minor logarithmic slowdown.

The budget-constrained wSAA estimator exhibits two distinct asymptotic behaviors for the following reasons.
With the budget allocation $m \sim \kappa \log\Gamma$ and $n \sim \Gamma/(\kappa\log\Gamma)$,
the statistical error in \eqref{eq:error-decomp} when scaled by $(\Gamma/(\kappa\log\Gamma))^{(1-\delta d_x)/2}$ converges to a normal distribution as stated in Theorem~\ref{thm:clt-cont}.
Additionally,
by \eqref{eq:opti-error-linear},
the optimization error when scaled by the same factor is of order
\begin{align*}
        \left(\frac{\Gamma}{\kappa \log\Gamma}\right)^{(1-\delta d_x)/2} \theta^{m}
        = \left(\frac{1}{\kappa \log\Gamma}\right)^{(1-\delta d_x)/2} \theta^{m - \kappa \log\Gamma} \Gamma^{(1-\delta d_x)/2-\kappa \log(1/\theta)},
\end{align*}
which vanishes only when $\kappa \ge \kappa^*$.
Hence, the CLT in \eqref{eq:thm-lin-i} is a joint outcome of the scaled statistical error (converging to normality) and the scaled optimization error (diminishing to zero).
When $\kappa < \kappa^*$, while the scaled statistical error maintains asymptotic normality, the scaled optimization error grows without bound. By choosing a smaller scaling factor $\Gamma^{\kappa \log(1/\theta)}$, we ensure the scaled statistical error converging to zero while keeping the scaled optimization error stochastically bounded.
This explains the asymptotic behavior in \eqref{eq:thm-lin-ii}.
However, characterizing the exact asymptotic distribution requires stronger assumptions than those used in our analysis.
Specifically, we would need additional requirements for the conditional distribution of $Y$, as well as more precise convergence guarantees for the optimization algorithm such as lower bounds on the per-iteration improvement.
Our current analysis relies only on upper bounds for the three convergent regimes under consideration.

\begin{figure}[ht]
\FIGURE{\includegraphics[width=\textwidth]{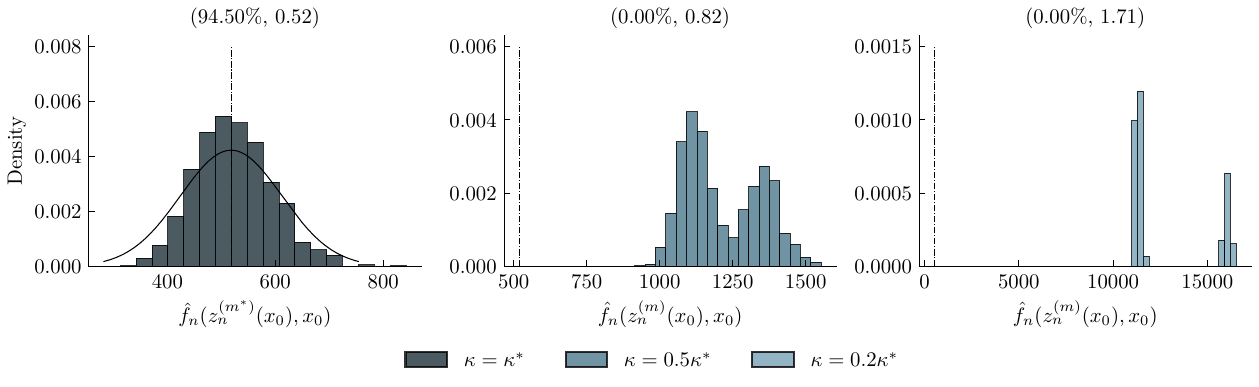}}
{Collapse of Asymptotic Normality of Budget-Constrained wSAA Estimator\label{fig:misspec}}
{Budget allocation is $m = \kappa \log \Gamma$.
The vertical dashed line represents $f^*(x_0)$.
The joint distribution of $(X,Y)$ and other relevant parameters are identical to those in Section~\ref{app:exp-fourthpoly} of the e-companion, except that the cost function $F(z; y)$ is a second-order polynomial rather than a fourth-order polynomial. The new covariate observation $x_0$ is selected to be at the 75\% quantile of the marginal distribution of $X$.
The computational budget is $\Gamma = 10^4$.
For each pair of parentheses, the first number represents the coverage of the asymptotic-normality-based 95\% confidence interval (Section~\ref{sec:ci}), while the second number indicates the relative width of the interval, normalized by $f^*(x_0)$.
}
\end{figure}

Figure~\ref{fig:misspec} illustrates, through numerical experiments, how the budget-constrained wSAA estimator $\hat{f}_n(z^{(m)}_n(x_0),x_0)$ exhibits distinct asymptotic behaviors depending on the value of $\kappa$. When $\kappa$ falls below the threshold $\kappa^*$, the optimization error dominates, possibly resulting in a non-normal asymptotic distribution that significantly complicates uncertainty quantification.

\section{
``Over-optimizing'' for  Normality
} \label{sec:overopt}

From Theorem \ref{thm:lin}, we know that for linearly convergent algorithms, the asymptotically optimal budget allocation follows $m(\Gamma) \sim \kappa^* \log\Gamma$. This allocation enables the budget-constrained wSAA estimator to attain a convergence rate nearly matching that of its unconstrained counterpart.

A practical challenge lies in determining $\kappa^*$, which depends on $\theta$---the convergence rate parameter of the optimization algorithm.
Example \ref{example:linear-conv} illustrates this challenge: for gradient descent methods applied to problems where $\hat{f}_n(\cdot, x_0)$ is $\lambda$-strongly convex with $L$-Lipschitz continuous derivatives, $\theta$ is determined by $\lambda$ and $L$. However, both parameters are typically difficult to estimate accurately.
The consequences of misspecifying $\theta$ can be severe. If $\theta$ is substantially underestimated (i.e., assumed to be close to zero), it falsely implies faster algorithmic convergence than is actually achieved. This misspecification leads to two  issues.
First, the budget-constrained wSAA estimator converges much slowly than the optimal rate of $(\Gamma/\log\Gamma)^{-(1 - \delta d_x)/2}$.
Second, and more critically, the error $\hat{f}_n(z^{(m)}_n(x_0),x_0) - f^*(x_0)$ may not be asymptotically normal anymore, invalidating asymptotic-normality-based uncertainty quantification.

To address potential misspecification of $\theta$, we propose an over-optimizing strategy for budget allocation, following the spirit of \cite{RoysetSzechtman12}.
The strategy uses an allocation rule of $m(\Gamma) \sim c_0 \Gamma^{\tilde{\kappa}}$ for some constants $c_0>0$ and $\tilde{\kappa}\in (0, 1)$.
Under this budget allocation, the number of iterations used to solve the wSAA problem is polynomial in $\Gamma$, which is significantly larger than the optimal quantity (i.e., $O(\Gamma^{\tilde{\kappa}})$ versus $O(\log\Gamma)$)---hence the concept of over-optimizing. The following theorem provides a formal statement of this strategy.

\begin{theorem} \label{thm:lin-overopt}
Suppose Assumptions~\ref{assump:regularity}--\ref{assump:uni-sol} hold, and $\scrA$ is linearly convergent  with parameter $\theta\in(0,1)$.
Consider an asymptotically admissible budget allocation $\{n(\Gamma),m(\Gamma)\}_{\Gamma \in \N_+}$ that satisfies $m(\Gamma) \sim c_0 \Gamma^{\tilde{\kappa}}$ as $\Gamma\to\infty$ for some constants $c_0 > 0$ and $\tilde{\kappa} \in (0,1)$.
Then,
    \begin{align*}
        \Gamma^{(1-\tilde{\kappa})(1-\delta d_x)/2} (\hat{f}_n(z^{(m)}_n(x_0),x_0) - f^*(x_0)) \Rightarrow \left(\frac{c_0^{1-\delta d_x}}{h_0^{d_x}}\right)^{1/2} N(0,\msfV(z^*(x_0),x_0)),
    \end{align*}
    as $\Gamma\to\infty$, for all $x_0 \in \scrX$,
    where $\msfV(z^*(x_0),x_0)$ is defined in \eqref{eq:V-expression}.
\end{theorem}

\begin{figure}[ht]
    \FIGURE{\includegraphics[width=\textwidth]{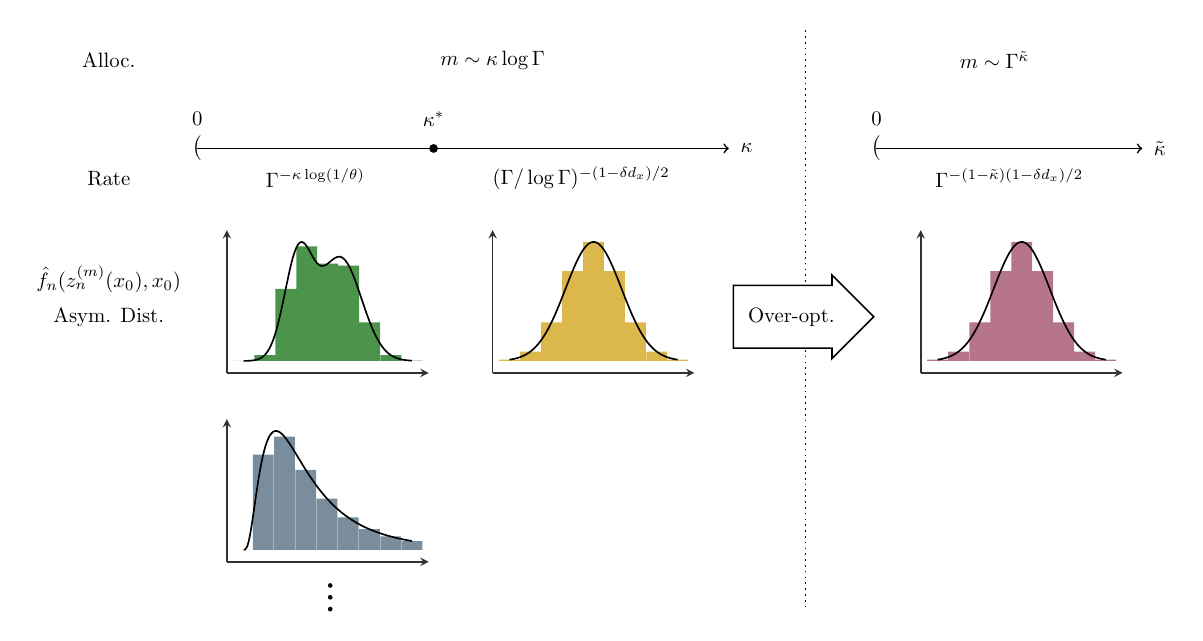}}
{Benefits of Over-optimizing \label{fig:overopt}}
{The optimization algorithm is linearly convergent. $\kappa^* = (1-\delta d_x)/(2\log (1/\theta))$ is the threshold value in Theorem~\ref{thm:lin} that determines three key properties: the optimal budget allocation, the optimal convergence rate of the budget-constrained wSAA estimator, and the conditions for asymptotic normality.}
\end{figure}

Over-optimizing offers two key benefits. First, the allocation rule does not depend on $\theta$, eliminating the need for its estimation. Second, it ensures asymptotic normality of the budget-constrained wSAA estimator regardless of the chosen $c_0$ and $\tilde{\kappa}$, thereby enabling more robust uncertainty quantification, see Figure~\ref{fig:overopt} for an illustration.

This strategy comes at a cost: when over-optimizing the wSAA problem, the resulting estimator's convergence rate $\Gamma^{(1-\tilde{\kappa})(1-\delta d_x)/2}$ is slower than the optimal rate $(\Gamma/\log\Gamma)^{-(1 - \delta d_x)/2}$. However, decision-makers can make this performance gap arbitrarily small by choosing a sufficiently small $\tilde{\kappa}$---a parameter under direct control---instead of relying on potentially inaccurate estimates of $\kappa^*$.
The essence of over-optimizing strategy is to accept a minimal, controllable reduction in convergence speed rather than risking an uncontrollable, potentially significant degradation in convergence rate caused by a misspecified $\kappa^*$, while simultaneously ensuring the validity of asymptotic-normality-based confidence intervals.

\section{Extensions} \label{sec:extension}

In this section, we extend the theory developed in Sections~\ref{sec:bdgt-alloc} and \ref{sec:overopt} to optimization algorithms of other convergent regimes: sublinearly and superlinearly convergent algorithms.

\subsection{Sublinearly Convergent Algorithms} \label{subsec:sublinear}

\begin{definition}[Sublinear Convergence] \label{defi:sublin}
    An algorithm $\scrA$ is said to converge sublinearly for solving the wSAA problem~\eqref{eq:cso-wSAA} if there exists a constant $\beta>0$ such that
    \begin{align}
        \hat{f}_n(z^{(m)}_n(x_0),x_0) - \hat{f}_n(\hat{z}_n(x_0),x_0) \le \frac{\Delta_n(x_0)}{m^\beta},  \label{eq:sublinear-conv}
    \end{align}
    for all $z_n^{(0)}(x_0) \in \scrZ$ and $n, m \in \N_+$,
    where $\Delta_n(x_0) \xrightarrow{p} \Delta(x_0)$ as $n \to \infty$ for some constant $\Delta(x_0) > 0$.
\end{definition}

\begin{example}
If $\hat{f}_n(\cdot, x_0)$ is convex and $L$-Lipschitz continuous, then (projected) subgradient descent with  a fixed stepsize of $\mathrm{diam}(\scrZ)/\sqrt{m+1}$
achieves sublinear convergence with $\beta=1/2$, where $\mathrm{diam(\scrZ)} = \sup_{z, z' \in \scrZ} \|z - z\| < \infty$ is the diameter of the compact set $\scrZ$ \citep[Theorem~3.2.2]{Nesterov18}.
Under Assumption~\ref{assump:regularity}-\ref{assump:F-lipschitz}, $\hat{f}_n(\cdot, x_0)$ satisfies this condition with $L=C_F$. In this case, $\Delta_n(x_0) = \mathrm{diam(\scrZ)} L$.
\end{example}

By matching the optimization error, which is bounded by \eqref{eq:sublinear-conv}, and the statistical error in the decomposition \eqref{eq:error-decomp}, an analysis analogous to that for linearly convergent algorithms implies that the asymptotically optimal budget allocation should satisfy $n^{-(1-\delta d_x)/2} \asymp m^{-\beta}$,  resulting in $m$ being a polynomial function of $\Gamma$. Theorem \ref{thm:sublin} characterizes the asymptotic behavior of the budget-constrained wSAA estimator under this polynomial allocation of $\Gamma$.

\begin{theorem} \label{thm:sublin}
Suppose Assumptions~\ref{assump:regularity}--\ref{assump:uni-sol} hold and $\scrA$ is sublinearly convergent with parameter $\beta>0$.
Consider an asymptotically admissible budget allocation $\{n(\Gamma),m(\Gamma)\}_{\Gamma \in \N_+}$ that satisfies $m(\Gamma) \sim c_0 \Gamma^{\kappa}$ as $\Gamma\to\infty$ for some constants $c_0 > 0$ and $\kappa \in (0,1)$.
    \begin{enumerate}[label=(\roman*)]
        \item
    If $(1-\delta d_x)/(1-\delta d_x+2\beta) \leq  \kappa < 1$, then
    \begin{align*}
        \Gamma^{(1-\kappa)(1-\delta d_x)/2}(\hat{f}_n(z^{(m)}_n(x_0),x_0) - f^*(x_0)) \Rightarrow \left(\frac{c_0^{1-\delta d_x}}{h_0^{d_x}}\right)^{1/2} N(0,\msfV(z^*(x_0),x_0)),
    \end{align*}
    as $\Gamma\to\infty$, for all $x_0 \in \scrX$,
    where $\msfV(z^*(x_0),x_0)$ is defined in \eqref{eq:V-expression}.
    \item If $0 < \kappa < (1-\delta d_x)/(1-\delta d_x+2\beta)$, then
    \begin{align*}
        \Gamma^{\kappa \beta} (\hat{f}_n(z^{(m)}_n(x_0),x_0) - f^*(x_0)) = O_{\pr}(1),
    \end{align*}
    as $\Gamma\to\infty$, for all $x_0 \in \scrX$.
    \end{enumerate}
\end{theorem}

Similar to the case of linear convergence, $\hat{f}_n(z^{(m)}_n(x_0),x_0)$ exhibits two types of asymptotic behavior depending on the budget allocation. Asymptotic normality arises when  $\kappa$ is large enough in the allocation $m(\Gamma) \sim c_0 \Gamma^\kappa$. The optimal allocation is to set $\kappa = (1-\delta d_x)/(1-\delta d_x+2\beta)$, yielding a convergence rate of $\Gamma^{-\beta(1-\delta d_x)/(1 - \delta d_x + 2\beta)}$.
In contrast to linearly convergent algorithms---where the optimal budget allocation almost recovers the unconstrained wSAA estimator's convergence rate of $n^{-(1-\delta d_x)/2}$, this rate is considerably slower. The computational budget constraint causes a more pronounced degradation in the convergence rate: a polynomial slowdown in $\Gamma$, as opposed to the logarithmic slowdown observed for linearly convergent algorithms. This stronger degradation results from the algorithm's slower convergence, which demands significantly greater computational effort to reduce the optimization error relative to the statistical error.

\begin{remark}
We do not explore over-optimizing strategy for sublinearly convergent algorithms since the optimal allocation parameters are more readily determined. Unlike the key parameter $\theta$ for linearly convergent algorithms, the key parameter $\beta$ in many sublinear algorithms can be clearly specified without estimating the hard-to-determine structural properties (e.g., Lipschitz constant) of the objective function. This mitigates the risk of significantly underestimating the threshold for optimal budget allocation, which could otherwise degrade the convergence rate of the budget-constrained wSAA estimator and break its asymptotic normality.
\end{remark}

\subsection{Superlinearly Convergent Algorithms} \label{subsec:superlinear}

\begin{definition}[Superlinear Convergence] \label{defi:superlin}
     An algorithm $\scrA$ is said to converge superlinearly for solving the wSAA problem~\eqref{eq:cso-wSAA} if there exist constants $\theta > 0$ and $\eta > 1$ such that
    \begin{align}
        \hat{f}_n(z^{(m)}_n(x_0),x_0) - \hat{f}_n(\hat{z}_n(x_0),x_0) \le \theta (\hat{f}_n(z^{(m-1)}_n(x_0),x_0) - \hat{f}_n(\hat{z}_n(x_0),x_0))^{\eta}, \label{eq:suplinear-conv}
    \end{align}
    for all $z^{(0)}_n(x_0) \in \scrZ$ and $n,m \in \N_+$.
\end{definition}

\begin{example}

If $\hat{f}_n(\cdot,x_0)$ is $\lambda$-strongly convex and twice differentiable with $L$-Lipschitz continuous second-order derivatives, then Newton's method achieves quadratic convergence ($\eta=2$), provided that the initial point is sufficiently close to the optimal solution $\hat{z}_n(x_0)$ \citep[Chapter 9.5]{Boyd04}. The projected Newton method under Hessian-induced norm
similarly achieves quadratic convergence \citep{Schmidt12}.
\end{example}

While the asymptotic analysis of the budget-constrained wSAA estimator follows similar principles across all convergent regimes, the case of superlinear convergence requires a more nuanced analysis and imposes additional restrictions on the algorithm's initial point.
By recursively applying inequality \eqref{eq:suplinear-conv}, we obtain the following bound on the optimization error:
\begin{align*}
    &\quad \hat{f}_n(z^{(m)}_n(x_0),x_0)  -
    \hat{f}_n(\hat{z}_n(x_0),x_0) \\
    &\leq \theta^{-1/(\eta-1)} \left( \theta^{1/(\eta-1)} ( \hat{f}_n(z^{(0)}_n(x_0),x_0) - \hat{f}_n(\hat{z}_n(x_0),x_0) ) \right)^{\eta^m} \\
    &= \theta^{-1/(\eta-1)} \Big[ \underbrace{\theta^{1/(\eta-1)} ( \hat{f}_n(z_n^{(0)}(x_0),x_0) - f(z_n^{(0)}(x_0),x_0))}_{\text{term 1}} + \underbrace{\theta^{1/(\eta-1)} ( f^*(x_0) - \hat{f}_n(\hat{z}_n(x_0),x_0) )}_{\text{term 2}} \\
    &\qquad + \underbrace{\theta^{1/(\eta-1)} ( f(z_n^{(0)}(x_0),x_0) - f^*(x_0)}_{\text{term 3}} ) \Big]^{\eta^m}.
\end{align*}
Here, both term 1 and term 2 are errors from approximating $f$ via $\hat{f}_n$, which can be quantified by Proposition~\ref{prop:fclt-cont} and Theorem~\ref{thm:clt-cont}, respectively.
Term 3, on the other hand, pertains to the distance between the initial solution $z_n^{(0)}(x_0)$ and the optimal solution $z^*(x_0)$ of the CSO problem \eqref{eq:cso}. To quantify the initial optimality gap, we define
\begin{align}
    \psi(z^{(0)}_n(x_0),x_0) \coloneqq \log(\theta^{-1/(\eta-1)} (f(z^{(0)}_n(x_0),x_0) - f^*(x_0))^{-1}), \label{eq:psi-def}
\end{align}
 where larger values of $\psi$ indicate closer proximity to the optimum.
Matching the optimization error and the statistical error in decomposition \eqref{eq:error-decomp} leads to the asymptotically optimal budget allocation---characterized by $n^{-(1-\delta d_x)/2} \asymp \exp(-\eta^m \psi(z^{(0)}_n(x_0),x_0))$, which implies $m \asymp \log\log\Gamma$. The following theorem formalizes this result.

\begin{theorem} \label{thm:superlin}
Suppose Assumptions~\ref{assump:regularity}--\ref{assump:uni-sol} hold, and $\scrA$ is superlinearly convergent with parameters $\theta > 0$ and $\eta > 1$. Further suppose $\psi(z^{(0)}_n(x_0),x_0) > 0$. Consider an asymptotically admissible budget allocation $\{n(\Gamma),m(\Gamma)\}_{\Gamma \in \N}$ that satisfies $m(\Gamma)\sim \kappa \log \log \Gamma$ as $\Gamma\to\infty$ for some constant $\kappa > 0$.
    \begin{enumerate}[label=(\roman*)]
        \item If either
        \begin{enumerate*}[label=(\alph*)]
        \item $\kappa > 1/\log \eta$, or
        \item $\kappa = 1/\log \eta$ and $\psi(z^{(0)}_n(x_0),x_0) \ge (1-\delta d_x)/2$,
        \end{enumerate*}
        then
        \begin{align*}
            \left(\frac{\Gamma}{\log\log\Gamma}\right)^{(1-\delta d_x)/2} (\hat{f}_n(z^{(m)}_n(x_0),x_0) - f^*(x_0)) \Rightarrow \left( \frac{\kappa^{1-\delta d_x}}{h_0^{d_x}} \right)^{1/2} N(0,\msfV(z^*(x_0),x_0)),
        \end{align*}
        as $\Gamma\to\infty$, for all $x_0 \in \scrX$,
        where $\msfV(z^*(x_0),x_0)$ is defined in \eqref{eq:V-expression}.
        \item If either
        \begin{enumerate*}[label=(\alph*)]
            \item $0< \kappa < 1/\log \eta$, or
            \item $\kappa = 1/\log \eta$ and $0 < \psi(z^{(0)}_n(x_0),x_0) < (1-\delta d_x)/2$,
        \end{enumerate*} then
        \begin{align*}
            \exp(\psi(z^{(0)}_n(x_0),x_0) (\log \Gamma)^{\kappa \log \eta})(\hat{f}_n(z^{(m)}_n(x_0),x_0) - f^*(x_0)) = O_{\pr}(1), %
        \end{align*}
        as $\Gamma\to\infty$, for all $x_0 \in \scrX$.
    \end{enumerate}
\end{theorem}

Theorem~\ref{thm:superlin} shows that for superlinearly convergent algorithms, the optimal convergence rate of the budget-constrained wSAA estimator is $(\Gamma/\log\log\Gamma)^{-(1-\delta d_x)/2}$, closely matching its unconstrained counterpart up to a double-logarithmic factor.
This optimal rate requires two conditions.
First, the budget allocation must satisfy $m(\Gamma) \sim \kappa^* \log\log\Gamma$, where $\kappa^* = 1/\log\eta$.
Second, the initial solution must be sufficiently close to the optimal solution $\hat{z}_n(x_0)$, with $\psi(z^{(0)}_n(x_0),x_0) \ge (1-\delta d_x)/2$.
However, if $\psi(z^{(0)}_n(x_0),x_0) < (1-\delta d_x)/2$,
then even under the asymptotically optimal budget allocation $m(\Gamma) \sim \kappa^* \log\log\Gamma$,
the convergence rate of $\hat{f}_n(z^{(m)}_n(x_0),x_0)$ deteriorates to
$\exp(-\psi(z^{(0)}_n(x_0),x_0) (\log \Gamma)^{\kappa^* \log \eta}) = \Gamma^{-\psi(z^{(0)}_n(x_0),x_0)}$,
which is much slower than the ideal rate of
$\Gamma^{-(1 - \delta d_x)/2}$.
Furthermore, it compromises the estimator's asymptotic normality, rendering asymptotic-normality-based uncertainty quantification invalid.

In practice, estimating $\psi = \psi(z^{(0)}_n(x_0),x_0)$ is challenging since it requires knowledge of the optimal conditional expected cost $f^*(x_0)$, precisely the quantity we seek to estimate in the first place.
To address this challenge, we apply the over-optimizing strategy developed in Section~\ref{sec:overopt} for linearly convergent algorithms, proactively sacrificing a negligible fraction of the convergence rate to prevent severe performance degradation and collapse of asymptotic normality.

\begin{theorem} \label{thm:superlin-overopt}
Suppose Assumptions~\ref{assump:regularity}--\ref{assump:uni-sol} hold, and $\scrA$ is superlinearly convergent with parameters $\theta > 0$ and $\eta > 1$. Further suppose $\psi(z^{(0)}_n(x_0),x_0) > 0$.
Consider an asymptotically admissible budget allocation $\{n(\Gamma),m(\Gamma)\}_{\Gamma \in \N_+}$ that satisfies $m(\Gamma) \sim \tilde{\kappa} \log \Gamma$ as $\Gamma\to\infty$ for some constant $ \tilde{\kappa} > 0$.
Then,
    \begin{align*}
        \left(\frac{\Gamma}{\log\Gamma}\right)^{(1-\delta d_x)/2} (\hat{f}_n(z^{(m)}_n(x_0),x_0) - f^*(x_0)) \Rightarrow \left( \frac{\tilde{\kappa}^{1-\delta d_x}}{h_0^{d_x}} \right)^{1/2} N(0,\msfV(z^*(x_0),x_0)),
    \end{align*}
    as $\Gamma\to\infty$, for all $x_0 \in \scrX$,
    where $\msfV(z^*(x_0),x_0)$ is defined in \eqref{eq:V-expression}.
\end{theorem}

Theorem \ref{thm:superlin-overopt} states that allocating the budget as $m(\Gamma) \sim \tilde{\kappa}\log \Gamma$ for some constant $\tilde{\kappa} > 0$, independent of the optimization algorithm's convergence parameters ($\eta$ and $\psi$), leads to two key results. First, the budget-constrained wSAA estimator converges at a rate of $(\Gamma/\log\Gamma)^{-(1-\delta d_x)/2}$, which is almost indistinguishable from the optimal rate of $(\Gamma/\log\log\Gamma)^{-(1-\delta d_x)/2}$. Second, the estimator preserves asymptotic normality, enabling more robust uncertainty quantification.

To conclude this section, Table~\ref{tab:summary} summarizes and compares our key findings on the budget-constrained wSAA estimator across  optimization algorithms of three convergent regimes.

\begin{table}[ht]
\TABLE
{A Summary of Theorems~\ref{thm:lin}--\ref{thm:superlin-overopt} \label{tab:summary}}
{
\begin{tabular}{ccccccccc}
    \toprule
    \multirow{3}{*}{$\scrA$} & \multicolumn{3}{c}{Optimal Allocation} & & \multicolumn{3}{c}{Over-optimizing} \\
    \cmidrule{2-4} \cmidrule{6-8}
    & $m^*$ & $\kappa^*$ & Rate & & $m$ & $\tilde{\kappa}$ & Rate \\
    \midrule
    Sublinear & $c_0 \Gamma^{\kappa^*}$ & $\displaystyle \frac{1-\delta d_x}{1-\delta d_x+2\beta}$ & $\Gamma^{-\kappa^*\beta}$ & & N/A & N/A & N/A \\[2ex]
    Linear & $\kappa^* \log \Gamma$ & $\displaystyle \frac{1-\delta d_x}{2\log (1/\theta)}$ & $\displaystyle \left(\frac{\Gamma}{\log\Gamma}\right)^{-(1-\delta d_x)/2}$ & & $c_0 \Gamma^{\tilde{\kappa}}$ & $(0,1)$ & $\Gamma^{-(1-\tilde{\kappa})(1-\delta d_x)/2}$ \\[2ex]
    Superlinear & $\kappa^* \log \log \Gamma$ & $\displaystyle \frac{1}{\log \eta}$ & $\displaystyle \left(\frac{\Gamma}{\log\log\Gamma}\right)^{-(1-\delta d_x)/2}$ &  & $\tilde{\kappa}  \log\Gamma$ & $(0,\infty)$ & $\displaystyle\left(\frac{\Gamma}{\log\Gamma}\right)^{-(1-\delta d_x)/2}$ \\[2ex]
    \bottomrule
\end{tabular}
}
{\emph{Note.} \begin{enumerate*}[label=(\roman*)]
    \item ``Rate'' means the rate of convergence in distribution of the budget-constrained wSAA estimator.
    \item $c_0 > 0$, $\theta \in (0,1)$, $\beta > 0$, and $\eta > 1$.
    \item ``N/A'' is the shorthand for ``not applicable''.
\end{enumerate*}}
\end{table}

\section{Confidence Intervals} \label{sec:ci}

The CLTs developed in Sections~\ref{sec:bdgt-alloc}--\ref{sec:extension} yield asymptotically valid confidence intervals for the optimal conditional expected cost $f^*(x_0)$. For example, Theorem~\ref{thm:clt-cont} implies that if the wSAA problem~\eqref{eq:cso-wSAA} can be solved to optimality, then for any $\alpha\in(0,1)$, an asymptotically valid $100(1-\alpha)\%$ confidence interval for $f^*(x_0)$ is given by
\begin{equation}\label{eq:ideal-CI}
\ci^\alpha(x_0) = \left[ \hat{f}_n(\hat{z}_n(x_0),x_0) \mp \Phi^{-1}(1-\alpha/2) \sqrt{\msfV (z^*(x_0), x_0)/(nh_n^{d_x})}\right].
\end{equation}
To implement this confidence interval, we need to estimate the typically unknown limiting variance $\msfV (z^*(x_0), x_0)$.
We accomplish this in two steps.
First, we construct a consistent estimator of $\msfV (z, x_0)$ for any given $z \in \scrZ$. Second, we plug in $z=\hat{z}_n(x_0)$ as a consistent estimator of $z^*(x_0)$, which yields a consistent estimator of $\msfV (z^*(x_0), x_0)$.

One could construct an estimator of $\msfV (z, x_0)$ directly from its definition in~\eqref{eq:V-expression}:
$\widetilde{\msfV}_n (z, x_0) \coloneqq  \hat{\sigma}^2_n (z,x_0) R_2(K)/\hat{p}_n(x_0)$,
where
\begin{equation} \label{eq:sigma_hat}
    \hat{\sigma}_n^2(z,x_0) \coloneqq  \sum_{i=1}^n w_{n}(x_i,x_0) (F(z;y_i) - \hat{f}_n(z,x_0))^2
\end{equation}
denotes the sample conditional variance and $\hat{p}_n(x_0) = (nh_n^{d_x})^{-1} \sum_{i=1}^n K ((x_i-x_0)/h_n)$
denotes the kernel density estimator of $p(x_0)$.
Although $\widetilde{\msfV}_n (z, x_0)$ can be shown to be a consistent estimator of $\msfV_n (z, x_0)$, it is often  numerically unstable due to the presence of $\hat{p}_n(x_0)$ in the denominator. This instability leads to large estimation variances, especially when only a few observations in $\scrD_n$ are near $x_0$.
Therefore, simply replacing $\msfV (z^*(x_0), x_0)$ in \eqref{eq:ideal-CI} with $\widetilde{\msfV}_n (\hat{z}_n(x_0), x_0)$ would result in excessively wide confidence intervals with substantial over-coverage.

To address this issue,
we consider the following estimator of $\msfV_n (z, x_0)$ without using $\hat{p}(x_0)$:
\begin{align} \label{eq:asymp-varest}
\widehat{\msfV}_n (z, x_0) \coloneqq (nh_n^{d_x}) \hat{\sigma}_n^2(z,x_0) \sum_{i=1}^n w_n^2(x_i,x_0).
\end{align}
It leads to the following confidence intervals for $f^*(x_0)$ based on the wSAA estimator $\hat{f}_n(\hat{z}_n(x_0),x_0)$ and its budget-constrained counterpart $\hat{f}_n(z^{(m)}_n(x_0),x_0)$, respectively:
\begin{align}
    \widehat{\ci}_{n}^\alpha(x_0) \coloneqq {}& \left[ \hat{f}_n(\hat{z}_n(x_0),x_0) \mp \Phi^{-1}(1-\alpha/2) \sqrt{\widehat{\msfV}_n (\hat{z}_n(x_0), x_0)/(nh_n^{d_x})}\right], \label{eq:ci-wSAA} \\
    \widehat{\ci}_{n,m}^\alpha(x_0) \coloneqq{}& \left[ \hat{f}_n(z^{(m)}_n(x_0),x_0) \mp \Phi^{-1}(1-\alpha/2) \sqrt{\widehat{\msfV}_n (z^{(m)}_n(x_0), x_0)/(nh_n^{d_x})}\right].\label{eq:ci-wSAA-constrained}
\end{align}
Corollaries~\ref{cor:ci-cvg} and \ref{cor:contrained-CI} show that they are asymptotically valid at the $100(1-\alpha)\%$ confidence level: $\lim_{n \to \infty} \pr( f^*(x_0) \in \widehat{\ci}_n^\alpha(x_0)) = 1-\alpha$ and $\lim_{\Gamma \to \infty} \pr( f^*(x_0) \in \widehat{\ci}_{n,m}^\alpha(x_0)) = 1-\alpha$.
The key to the proofs lies in showing the consistency of $\hat{\sigma}^2_n (\hat{z}_n(x_0),x_0)$ as an estimator of $\sigma^2 (z^*(x_0),x_0)$, then using the fact that $nh_n^{d_x} \sum_{i=1}^n w_n^2(x_i,x_0) = R_2(K)/p(x_0) + o_{\pr}(1)$.

\begin{corollary} \label{cor:ci-cvg}
    Consider the setting of the CLT in Theorem~\ref{thm:clt-cont}. Suppose
    $\int_{\RR^{d_x}} K^4(u) du < \infty$, $\E[|F^2(z;Y)| (\log |F^2(z;Y)|)^{+}] < \infty$ for all $z \in \scrZ$, and the wSAA problem~\eqref{eq:cso-wSAA} has a unique solution $\hat{z}_n(x_0)$.
    Then,
    \begin{align*}
        \frac{\hat{f}_n(\hat{z}_n(x_0),x_0) - f^*(x_0)}{\sqrt{\hat{\sigma}_n^2(\hat{z}_n(x_0),x_0) \sum_{i=1}^n w_n^2(x_i,x_0)}} \Rightarrow N(0,1),
    \end{align*}
    as $n\to\infty$, for  all $x_0 \in \scrX$, where $\hat{\sigma}_n^2(z,x_0)$ is defined in \eqref{eq:sigma_hat}.
\end{corollary}

\begin{corollary}\label{cor:contrained-CI}
    Consider the setting of the CLT in any of Theorems~\ref{thm:lin}, \ref{thm:lin-overopt}, \ref{thm:sublin}, \ref{thm:superlin}, or \ref{thm:superlin-overopt}.
    Suppose $\int_{\RR^{d_x}} K^4(u) du < \infty$,
    $\E[|F^2(z;Y)| (\log |F^2(z;Y)|)^{+}] < \infty$ for all $z \in \scrZ$, and  the wSAA problem~\eqref{eq:cso-wSAA} has a unique solution $\hat{z}_n(x_0)$.
    Then,
    \begin{align*}
        \frac{\hat{f}_n(z^{(m)}_n(x_0), x_0) - f^*(x_0)}{\sqrt{\hat{\sigma}_n^2(z^{(m)}_n(x_0),x_0) \sum_{i=1}^n w_n^2(x_i,x_0)}} \Rightarrow N(0,1),
    \end{align*}
    as $\Gamma\to\infty$, for all $x_0 \in \scrX$, where $\hat{\sigma}_n^2(z,x_0)$ is defined in \eqref{eq:sigma_hat}.
\end{corollary}

\section{Case Study: Capacity Management for Bike Sharing} \label{sec:exper}

In this section, we conduct a numerical experiment using real data from a large bike-sharing platform to validate our theoretical results on a capacity management problem. We demonstrate the construction of confidence intervals for the optimal conditional expected cost in two scenarios: the idealized case without computational constraints and the case where the wSAA problem is solved using a linearly convergent algorithm under computational budget. We also illustrate the benefits of over-optimizing strategy in the latter scenario. Additional numerical experiments for other CSO problems, solved with sublinearly and superlinearly convergent algorithms, are available in Sections~\ref{app:exp-nv} and \ref{app:exp-fourthpoly} of the e-companion, respectively.
The experimental design is detailed in Section~\ref{app:exp-design} of the e-companion.
All experiments are implemented in Python on a computer with a 3.1GHz AMD CPU and 32GB of RAM.
When needed, the solver Gurobi (\url{https://www.gurobi.com/}) is used to compute $f^*(x_0)$ and $\hat{f}_n(\hat{z}_n(x_0),x_0)$.

\subsection{Problem Description}

\begin{figure}[ht]
    \FIGURE{\includegraphics[width=230pt]{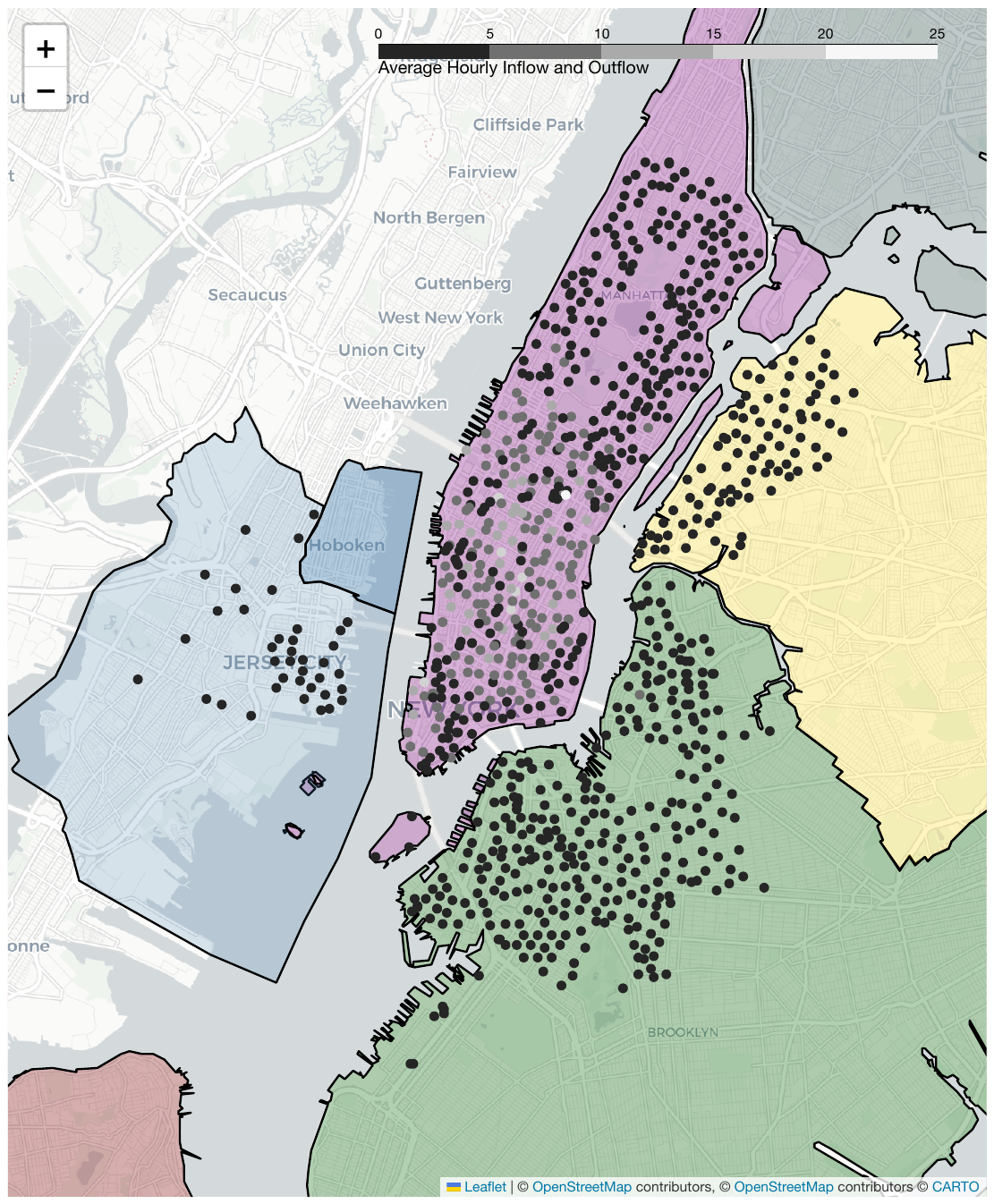}}
{Citi Bike Sharing System \label{fig:station}}
{Shaded regions represent the five boroughs of NYC---Manhattan (pink), Brooklyn (green), Queens (yellow), Bronx (gray), and Staten Island (red)---as well as Hoboken (darker blue) and Jersey City (lighter blue) in New Jersey. The hourly average inflow and outflow at each station is calculated over the horizon from June 1, 2013 to December 31, 2019,  after basic data cleaning.}
\end{figure}

Launched in May 2013, Citi Bike is the largest bike-sharing system in the United States, initially operating in New York City (NYC) and later expanding into New Jersey. Station locations are marked by dots in Figure~\ref{fig:station}, with color intensity encoding the magnitude of average inflow and outflow. Throughout our study period, a total of 926 stations were active, distributed across Manhattan (489), Brooklyn (312), Queens (87), and Jersey City (38). As shown in Figure~\ref{fig:station}, bike usage is uneven across the city, with most trips concentrated in Manhattan and Brooklyn, while bikes remain underutilized in other boroughs. This imbalance causes some stations being overcrowded while others are largely empty. In this case study, we investigate a critical capacity management problem: determining the optimal number of bikes allocated to each borough (e.g., Manhattan) during rush hours.

For hourly demand $Y$ and capacity $z$ (representing the initial number of bikes needed to accommodate uncertain rental requests) in a given borough, we consider the following cost function:
\begin{align}\label{eq:piecewise-quad}
F(z;Y) = c_u\left[(Y-z)^2 \ind(Y \geq z)\right] + c_o\left[(z-Y)^2 \ind(Y < z) \right],
\end{align}
where $c_u$ and $c_o$ denote the per-unit penalties for shortages (i.e., empty stations preventing rentals) and congestion (i.e., full stations preventing returns), respectively.
This cost function, adapted from \cite{Donti17}, penalizes shortages and congestion quadratically, making it more stringent than the standard newsvendor-type cost function. A modest surplus of bikes provides a buffer against unexpected demands, thereby enhancing operational flexibility. However, a significant surplus increases maintenance and depreciation costs while also raising congestion risks. Likewise, while small shortages may cause minor inconvenience as users can typically find bikes at nearby stations, substantial shortages can lead to long waiting times, thus frustrating commuters who rely on  timely services. Frequent bike unavailability may drive users toward alternative modes of transportation. This not only reduces short-term revenue but also erodes long-term market share due to increased customer churn.
The cost function~\eqref{eq:piecewise-quad} has a strong connection to expectile regression, which is a variant of quantile regression with penalty in the form of $\ell_1$ norm replaced by $\ell_2$ norm \citep{NeweyPowell87}. The CSO problem \eqref{eq:cso} can be solved explicitly, with the optimal solution equal to the conditional expectile at level $c_u/(c_u+c_o)$ of $Y$ given $X=x_0$. The optimal solution to the corresponding wSAA problem~\eqref{eq:cso-wSAA} has a similar interpretation.

\subsection{Data and CWGAN Simulator}

We constructed hourly demand data $Y$ using trip records from the Citi Bike website (\url{https://citibikenyc.com/system-data}) for Manhattan, covering the period from June 1, 2013 to December 31, 2019.
For simplicity, we treat the observed number of rentals as the true demand, though this ignores potential right-censoring when demand exceeds bike availability.
Our data cleaning process proceeds with pre-filtered records that exclude trips less than one minute (indicating equipment issues) or more than 48 hours (considered as lost/stolen per rental agreement).
We then apply additional filters, keeping only trips that start or end in Manhattan while excluding those occurring on holidays, weekends, during winter months, or under extreme weather conditions. We focus on ten busy hours: 8:00-10:00 AM and 12:00-8:00 PM. The final dataset contains 12,502 observations.

\begin{figure}[ht]
    \FIGURE{\includegraphics[width=0.5\textwidth]{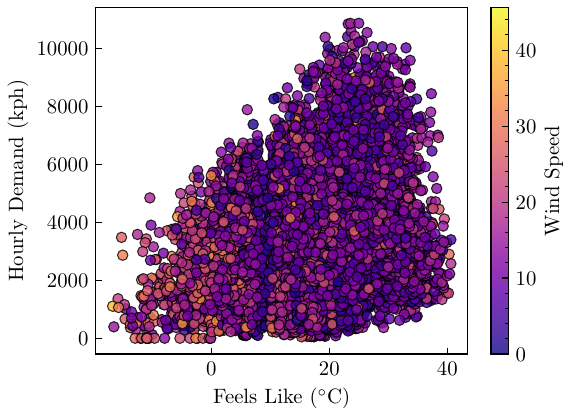}}
{``Feels Like'' Temperature, Wind Speed, and Hourly Demand \label{fig:citibike-data-cleaned}}
{}
\end{figure}

The covariate data $X$ is constructed using weather conditions obtained from Visual Crossing (\url{https://www.visualcrossing.com/}), a leading weather data provider.
We synchronize weather measurements with hourly demand data for temporal alignment.
While Visual Crossing offers comprehensive measurements including temperature, precipitation, humidity, wind speed/direction, and cloud cover, we found that only two of them demonstrate sufficient explanatory power for hourly demand variations.
The two measurements are apparent temperature (also known as ``feels like'' temperature) and wind speed, which together form a two-dimensional covariate vector (i.e., $d_x=2$).
Figure~\ref{fig:citibike-data-cleaned} illustrates their relationship with hourly demand:
demand tends to peak when the ``feels like'' temperature is between $10^\circ \text{C}$ and $30^\circ \text{C}$, and the wind speed is low.

\begin{figure}[ht]
\FIGURE{\includegraphics[width=320pt]{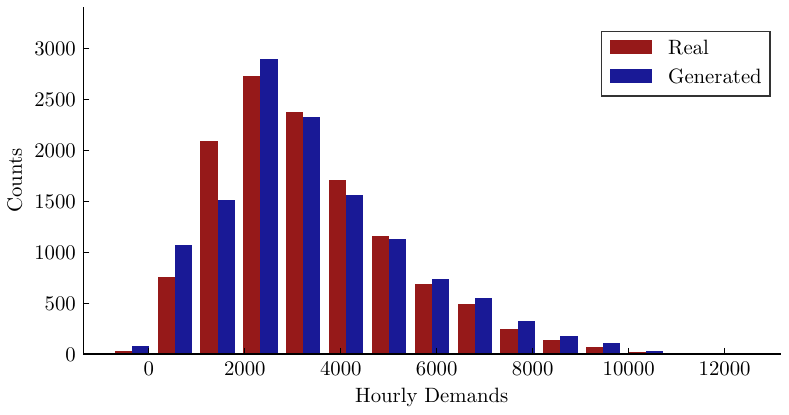}}
{Marginal Histograms of Hourly Demand \label{fig:citibike-CWGAN-perf}}
{}
\end{figure}

Since the true conditional distribution of $Y$ given $X$ is unknown, we adopt the approach from \cite{Athey24} to train a Conditional Wasserstein Generative Adversarial Network (CWGAN) on our dataset (see Section~\ref{app:CWGAN} of the e-companion for details). Its capability to reproduce the true demand is validated through the marginal distributions shown in Figure~\ref{fig:citibike-CWGAN-perf}.
The trained CWGAN then serves as a conditional distribution simulator, with its outputs treated as the ground truth. We then generate samples from this simulator to construct the dataset $\scrD_n$ in each macro-replication.

\subsection{Results}

We set the cost parameters to $c_u=1$ and $c_o=0.5$, and let the new covariate observation be $x_0 = (13.20^\circ \text{C}, 8.40 \text{kph})^\intercal$,
corresponding to the 75\% quantile ($\tau=0.75$) of the empirical marginal distributions of ``feels like'' temperature and windspeed, respectively.
We additionally explored  alternative parameter configurations, specifically $\tau \in \{0.25, 0.50\}$ and $c_o/c_u \in \{0.2, 0.8\}$. The experimental results are consistent across different configurations;
details are omitted due to space limit. In the wSAA problem~\eqref{eq:cso-wSAA},
we use the Gaussian kernel function
$K(u) = \exp(-\|u\|^2/2)$.
For the bandwidth $h_n = h_0 n^{-\delta}$, we set $\delta = 1/(d_x+3)=1/5$ and select $h_0$ via cross-validation (see Section~\ref{app:CV} of the e-companion). We set the confidence level to 95\% ($1-\alpha = 0.95$).

The first goal of this experiment is to evaluate the performance of the confidence interval~\eqref{eq:ci-wSAA} in terms of its width and coverage in the idealized case.
To obtain $\hat{f}_n(\hat{z}_n(x_0), x_0)$, the wSAA problem~\eqref{eq:cso-wSAA} is solved to optimality.
To compute the true optimal value $f^*(x_0)$ of the CSO problem~\eqref{eq:cso}, we generate a large number ($10^7$) of samples of $Y$ given $X=x_0$ from the CWGAN simulator, and then solve the resulting SAA problem.
For each sample size $n$, we compute three quantities based on 1,000 replications:
\begin{enumerate}[label=(\roman*)]
    \item relative width of the confidence interval (ratio of interval width to $f^*(x_0)$);
    \item coverage (frequency of $f^*(x_0)$ falling within the confidence interval); and
    \item relative root mean squared error (RMSE) (RMSE of $\hat{f}_n(\hat{z}_n(x_0), x_0)$ normalized by $f^*(x_0)$).
\end{enumerate}

\begin{figure}[ht]
    \FIGURE{\includegraphics[width=\textwidth]{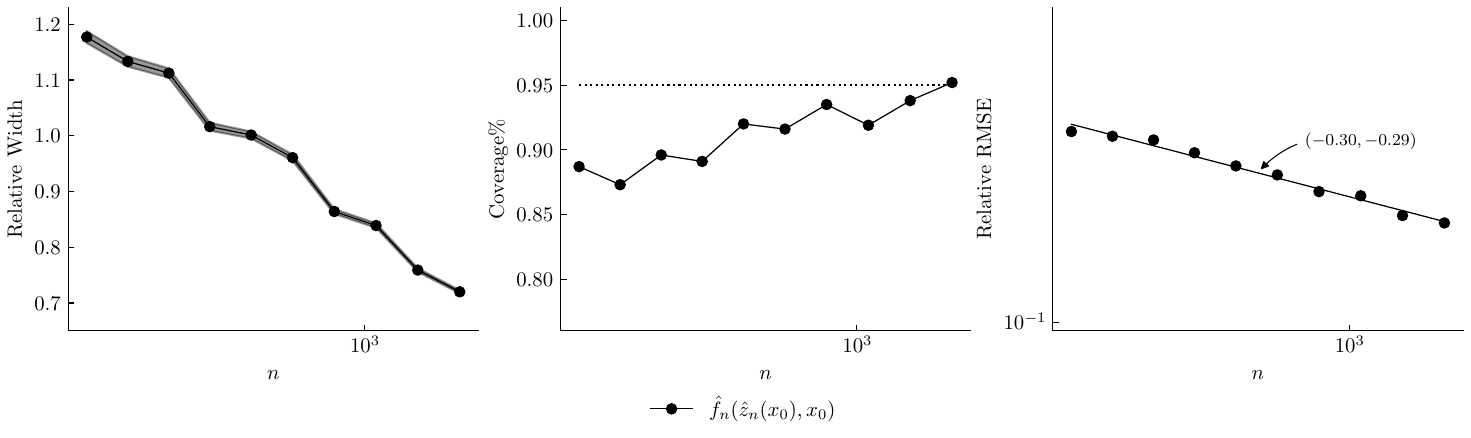}}
{Bike Sharing Problem (Without Computational Budget Constraint) \label{fig:exp-bike-o_0.5-tau_0.75-ideal}}
{The first number in the annotated parentheses represents the theoretical convergence rate of relative RMSE, which is $n^{-(1-\delta d_x)/2}$ (Theorem~\ref{thm:clt-cont}), while the second number indicates the empirical slope obtained from regressing log relative RMSEs on $\log n$.
}
\end{figure}

As shown in Figure~\ref{fig:exp-bike-o_0.5-tau_0.75-ideal}, when solving the wSAA problem~\eqref{eq:cso-wSAA} to optimality, the confidence interval~\eqref{eq:ci-wSAA} has asymptotically exact coverage of $f^*(x_0)$. With sample sizes around $10^3$, the empirical coverage approaches the target 95\% level (indicated by horizontal dashed line). As $n$ increases, the relative interval widths decrease with diminishing variability, appearing as a narrowing band---which represents the mean plus or minus one standard error---in the left panel. Moreover, the relative RMSEs of the optimal value estimates converge at a rate of $n^{-0.29}$, consistent with the theoretical rate of $n^{-0.30}$.

\begin{figure}[ht]    \FIGURE{\includegraphics[width=\textwidth]{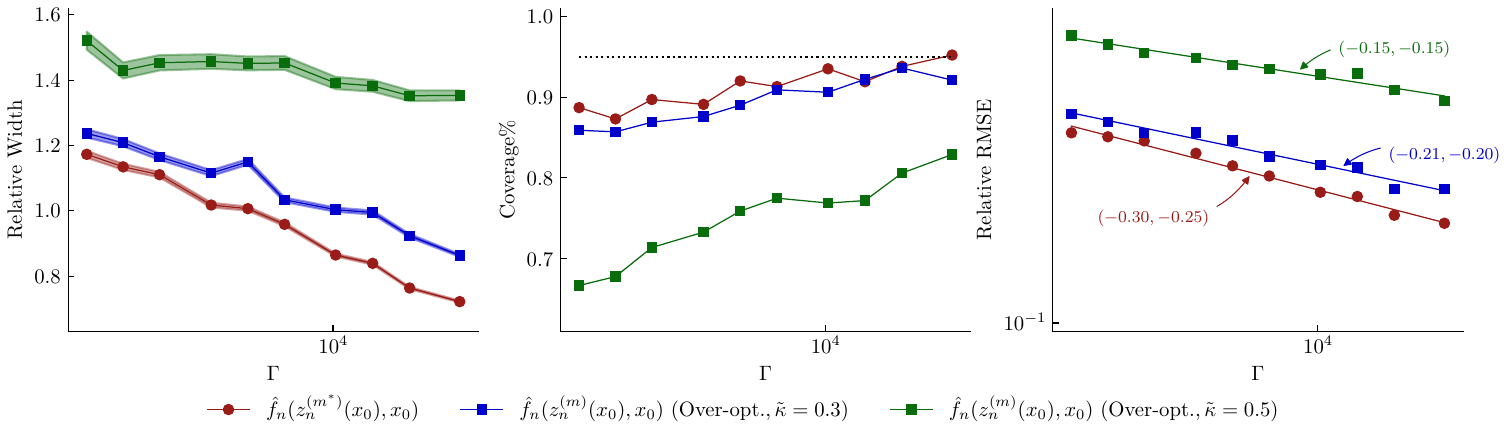}}
{Bike Sharing Problem (With Computational Budget Constraint)\label{fig:exp-bike-o_0.5-tau_0.75-bdgt}}
{The first number in each pair of annotated parentheses represents the theoretical convergence rate of relative RMSE,
while the second number indicates the empirical slope obtained from regressing log relative RMSEs on $\log \Gamma$.
Under the optimal budget allocation ($m = \kappa^* \log\Gamma$), the theoretical rate is $\Gamma^{-(1-\delta d_x)/2}$ (up to a logarithmic factor) (Theorem~\ref{thm:lin}). For the over-optimizing strategy ($m^* = \kappa^* \Gamma^{\tilde{\kappa}}$ with $\tilde{\kappa} = 0.3, 0.5 $), the rate is $\Gamma^{-(1-{\tilde{\kappa}})(1-\delta d_x)/2}$ (Theorem~\ref{thm:lin-overopt}). The optimization algorithm is linearly convergent.
}
\end{figure}

Next, we evaluate the performance of the confidence interval~\eqref{eq:ci-wSAA-constrained} when the wSAA problem~\eqref{eq:cso-wSAA} is solved under computational budget. Here, $\scrA$ denotes projected gradient descent with backtracking, where the line‑search parameters are set to $a=0.45$ and $b=0.9$.
The cost function~\eqref{eq:piecewise-quad} induces a $\lambda$-strongly convex objective $\hat{f}_n(\cdot,x_0)$ with $L$-Lipschitz continuous derivatives, where $\lambda = 2c_o$ and $L = 2c_u$. Therefore, $\scrA$ converges linearly when solving the wSAA problem (see Section~\ref{app:alg-para} of the e-companion for details).
For each computational budget $\Gamma$, we compare three allocation rules:
\begin{enumerate}[label=(\roman*)]
    \item $m(\Gamma) \sim \kappa^* \log\Gamma$ with $\kappa^* =  (1-\delta d_x)/(2\log(1/\theta))$;
    \item $m(\Gamma) \sim \kappa^* \Gamma^{\tilde{\kappa}}$ with $\tilde{\kappa}=0.3$; and
    \item $m(\Gamma) \sim \kappa^* \Gamma^{\tilde{\kappa}}$ with $\tilde{\kappa}=0.5$.
\end{enumerate}
The first rule is theoretically optimal, whereas the latter two represent over-optimizing strategies.
For each allocation rule, we then evaluate the relative width and coverage of the confidence interval~\eqref{eq:ci-wSAA-constrained}, as well as the relative RMSE of $\hat{f}_n(z^{(m)}_n(x_0),x_0)$.
The results are presented in Figure~\ref{fig:exp-bike-o_0.5-tau_0.75-bdgt}.

Under optimal budget allocation (represented by the red lines in Figure~\ref{fig:exp-bike-o_0.5-tau_0.75-bdgt}), the wSAA estimator performs nearly as well as it does without computational constraints. The confidence intervals exhibit a pattern---both in terms of width and coverage---similar to that observed in the idealized case (see Figure~\ref{fig:exp-bike-o_0.5-tau_0.75-ideal}).
The relative RMSE of $\hat{f}_n(z^{(m)}_n(x_0),x_0)$ converges at a rate of approximately $\Gamma^{-0.25}$, which is slightly slower than the theoretical rate of $\Gamma^{-0.30}$ (up to a logarithmic factor).

The over-optimizing strategy is evaluated as follows.
A budget allocation with $m \sim \kappa^* \Gamma^{0.5}$ (represented by the green lines in Figure~\ref{fig:exp-bike-o_0.5-tau_0.75-bdgt}) hurts the wSAA estimator's performance. Specifically, the confidence intervals are roughly 0.5 times wider than those under the optimal budget allocation, and the variations in their widths show no clear decreasing trend even when the budget exceeds $10^4$. Their coverage declines to around 80\%.
Moreover, the relative RMSEs are much larger and converge at a slower rate of $\Gamma^{-0.15}$, compared with the more favorable rate of $\Gamma^{-0.25}$. These observations indicate a poor statistical--computational tradeoff:
insufficient computing resources are allocated to estimating the objective function, resulting in inaccurate estimates, whose errors are further magnified in the subsequent optimization procedure.
However, as the over-optimizing strategy approaches the optimal allocation rule---specifically, when $\tilde{\kappa}$ is reduced from 0.5 to 0.3---its performance gap relative to the optimal allocation narrows.
Notably, the coverage of the confidence intervals becomes comparable to that of the optimal allocation once $\Gamma$ exceeds $10^4$. This finding demonstrates that a moderate degree of over-optimizing the problem, while initially suboptimal, can ultimately achieve a level of efficiency comparable to the optimal allocation as the computational budget increases.
Moreover, when the parameter $\theta$---which governs the algorithm's convergence behavior and thereby determines the budget allocation---cannot be precisely estimated, the over-optimizing strategy still guarantees the validity of uncertainty quantification and achieves a satisfactory convergence rate of $\Gamma^{-0.20}$.

\section{Conclusions} \label{sec:concl}

In this paper, we studied uncertainty quantification for CSO under computational constraints, focusing on the widely used wSAA method. We established CLTs for the budget-constrained wSAA estimator and constructed confidence intervals for optimal conditional expected costs that account for limited computational resources at solve time. We characterized the statistical--computational tradeoff faced by practitioners with historical data but finite budgets, and showed how to allocate resources between sample size (to reduce estimation error) and iterations (to reduce optimization error) based on the optimization algorithm’s convergence behavior.

Because the optimal allocation rule relies on structural parameters of the objective that may be misspecified and can break the asymptotic normality, we showed that modest ``over-optimization'' (a slight increase in iterations over the nominal allocation) can improve the reliability of uncertainty quantification without causing uncontrolled slowdowns.

Despite rapid progress in CSO \citep{Sadana25}, both uncertainty quantification and explicit treatment of computational constraints remain largely underexplored.
While this paper contributes to narrowing the gap, several promising research directions remain. First, extending our analysis beyond i.i.d. data to more general data-generating processes, such as Markov chains, would broaden its scope of applicability.
Second, our analysis could be extended to two-stage CSO problems where wSAA has been applied \citep{NotzPibernik22,BertsimasMcCordSturt23}.
Lastly, developing uncertainty quantification procedures for alternative CSO methods---particularly those that learn the policies $z(x_0)$ directly
\citep{BertsimasKoduri22}, as opposed to following the estimate-then-optimize paradigm of wSAA---constitutes another important direction.

\begin{APPENDICES}{}  \label{apx:tech-detail}

\end{APPENDICES}

\bibliographystyle{informs2014} %
\bibliography{ref} %

\ECSwitch

\EquationsNumberedBySection

\ECHead{Supplemental Material}

\section{Omitted Proofs} \label{ec:proof}

\subsection{Proof of Theorem~\ref{thm:clt-cont}}
Define $r(z,x_0) \coloneqq \int_{\scrY} F(z;y) p(x_0,y) d y$,
where $p(x,y)$ denotes the joint density of $X$ and $Y$. Then, we can write $f(z,x_0) = r(z,x_0)/p(x_0)$.

\subsubsection{A Technical Lemma}

\begin{lemma} \label{lem:fclt-r}
    $\sqrt{nh_n^{d_x}} (\hat{r}_n(\cdot,x_0)-r(\cdot,x_0))$ converges to a Gaussian process in distribution as $n\to\infty$.
\end{lemma}

\proof{Proof of Lemma~\ref{lem:fclt-r}.}
    Define $q_{ni}(z) \coloneqq \frac{1}{\sqrt{nh_{n}^{d_x}}} K\left(\frac{X_i-x_0}{h_n}\right) F(z;Y_i)$.
    We have $\hat{r}_n(z,x_0) = \sum_{i=1}^n q_{ni}(z)$, see Section~\ref{subsec:tradeoff-ideal}.
    The expectation of $q_{ni}(z)$ can be calculated as
    \begin{align} \label{eq:pf-fclt-exp-qni}
        \E[q_{ni}(z)]
        &= \sqrt{\frac{h_n^{d_x}}{n}} \int_{\RR^{d_x}} \frac{1}{h_n^{d_x}} K\left(\frac{x_i-x_0}{h_n}\right) f(z,x_i) p(x_i) d x_i = \sqrt{\frac{h_n^{d_x}}{n}} \int_{\RR^{d_x}} K(u) r(z,x_0+h_n u) d u \nonumber \\
        &= \sqrt{\frac{h_n^{d_x}}{n}} \int_{\RR^{d_x}} K(u) \left( r(z,x_0)+ h_n u^\intercal \nabla_x r(z,x_0)+ \frac{h_n^2}{2} u^\intercal \nabla_x^2 r(z,\bar{x}_0) u \right) d u \nonumber \\
        &= \sqrt{\frac{h_n^{d_x}}{n}}\left(r(z,x_0)+ \frac{h_n^2}{2} \int_{\RR^{d_x}} K(u) u^\intercal \nabla_x^2 r(z,\bar{x}_0) u d u \right),
    \end{align}
    where $\bar{x}_0$ is on the line segment joining $x_0$ and $x_0+h_n u$.
    The differentiability of $r(z,x)$ in $x$ for all $z \in \scrZ$ with bounded first- and second-order derivatives can be deduced from
    Assumption~\ref{assump:smooth-1}.

    We decompose the partial-sum process $\hat{r}_n(z, x_0)$ as follows:
    \begin{align*}
        \sqrt{nh_n^{d_x}} \left(\hat{r}_n(z, x_0)-r(z,x_0)\right)
        &= \underbrace{\left(\sum_{i=1}^n \left(q_{ni}(z) - \E[q_{ni}(z)]\right) \right)}_{\coloneqq I_{n1}(z,x_0)} + \underbrace{\left(\sum_{i=1}^n \E[q_{ni}(z)]-\sqrt{nh_n^{d_x}} r(z,x_0)\right)}_{\coloneqq I_{n2}(z,x_0)}.
    \end{align*}
    First, we can show that $\sup_{z \in \scrZ}|I_{n2}(z,x_0)| = (\sqrt{n h_n^{d_x+4}}/2) \sup_{z \in \scrZ}\left|\int_{\RR^{d_x}} K(u) u^\intercal \nabla_x^2 r(z,\bar{x}_0) u d u\right| \leq 2C_f C_p \Upsilon(K)\sqrt{n h_n^{d_x+4}} = C_1 \sqrt{n h_n^{d_x+4}} = o(1)$,
    where $C_1>0$ is a constant.
    The inequality is by Assumption~\ref{assump:smooth-1}, since the norm of $\nabla_x^2 r(z,\bar{x}_0) = \nabla_x^2 f(z,\bar{x}_0) p(\bar{x}_0) + f(z,\bar{x}_0) \nabla_x^2 p(\bar{x}_0) + 2 \nabla_x f(z,\bar{x}_0) \nabla_x p(\bar{x}_0)$ is bounded by $4C_f C_p$.
    Next, we argue that $I_{n1}(\cdot,x_0)$ converges to a Gaussian process in distribution by verifying the five conditions in Theorem 10.6 of \cite{Pollard90_ec}.

    \textbf{Verifying Condition (i).} %
        Let $\scrF_{n\omega} = \{F(z;Y_{ni}): z \in \scrZ, n \ge 1, 1 \le i \le n\}$ represent the collection of processes associated with the random function $F$, where $\omega$ represents the realizations of $\scrD_n$; that is, $((X_{n1}, Y_{n1}), \ldots, (X_{nn}, Y_{nn}))$. The inclusion of $n$ in the subscript is for constructing a triangular array. The envelope of $\scrF_{n\omega}$, denoted by $\{H_n(\omega) = (H(Y_{n1}), \ldots, H(Y_{nn})), n \ge 1\}$, can be constructed by
        $ H(Y_{ni}) = \max\left\{ \sup_{z \in \scrZ}|F(z;Y_{ni})|, C_F \right\} = \max\left\{M(Y_{ni}), C_F\right\}$.
        By Assumptions~\ref{assump:smooth-1} and \ref{assump:envelope}, $H(Y_{ni})$ is finite.
        For any nonnegative vector $\alpha = (\alpha_1, \ldots, \alpha_n)^\intercal \in \RR_+^n$, it holds that $\left\|\left(\alpha_1 F(z;Y_{n1}), \ldots, \alpha_n F(z;Y_{nn})\right) - \left(\alpha_1 F(z';Y_{n1}), \ldots, \alpha_n F(z';Y_{nn})\right)\right\| \le \sum_{i=1}^n \alpha_i C_F \|z -z'\| \le \sum_{i=1}^n \alpha_i H(Y_{ni}) \|z -z'\| = \|\alpha \odot H_n(\omega)\| \|z-z'\|$, where $\odot$ denotes the element-wise product.
        By Definition~3.3 of \cite{Pollard90_ec},
        $D(\epsilon \|\alpha \odot  H_n(\omega)\|, \alpha \odot \scrF_{n\omega}) \le D(\epsilon, \scrZ)$
        for any $\epsilon \in (0,1)$, where $D$ denotes the $\ell_2$ packing number.
        Since $\scrZ$ is compact, there exist constants $\widetilde{C}_1$ and $\widetilde{C}_2$ such that $D(\epsilon, \scrZ) \le \widetilde{C}_1\epsilon^{-\widetilde{C}_2}$ \citep[pp.18--20]{Pollard90_ec}. It then follows that $\int_0^1 \sqrt{\log D(\epsilon, \scrZ)} d\epsilon \le \sqrt{\log \widetilde{C}_1} +  \sqrt{\widetilde{C}_2} \int_0^{+\infty} t^{1/2} e^{-t} dt = \sqrt{\log \widetilde{C}_1} + \sqrt{\widetilde{C}_2 \pi}/2 < \infty$.
        Therefore, $\scrF_{n\omega}$
        is manageable in the sense of Definition~7.9 of \cite{Pollard90_ec}.
        Let $\{Q_{ni}: n \ge 1, 1 \le i \le n\}$ be the array of envelopes defined as $Q_{ni} \coloneqq \frac{1}{\sqrt{nh_{n}^{d_x}}} K\left(\frac{X_i-x_0}{h_n}\right) H(Y_i)$.
        Then, it can be easily verified that the random processes $\{q_{ni}(z): z \in \scrZ, n \ge 1, 1 \le i \le n\}$ are manageable with respect to $\{Q_{ni}\}$.

    \textbf{Verifying Condition (ii).} For any $z,z'\in \scrZ$, it holds that
        \begin{align*}
            &\quad \sum_{i=1}^n \E[q_{ni}(z) q_{ni}(z')]
            = \int_{\RR^{d_x}} \frac{1}{h_n^{d_x}} K^2\left(\frac{x_i-x_0}{h_n}\right) \nu(z,z',x_i) p(x_i) d x_i \\
            &= \int_{\RR^{d_x}} K^2(u) \nu(z,z',x_0+h_n u) p(x_0+h_n u) d u \\
            &= \int_{\RR^{d_x}} K^2(u) \nu(z,z',x_0) p(x_0) d u + h_n \int_{\RR^{d_x}} K^2(u) u^\intercal \left( \nabla_x \nu(z,z',\bar{x}_0) p(\bar{x}_0) + \nu(z,z',\bar{x}_0) \nabla_x p(\bar{x}_0) \right) d u \\
            &= R_2(K) \nu(z,z',x_0) p(x_0) + o(1).
        \end{align*}
        In addition, $\sum_{i=1}^n \E[q_{ni}(z)] \E[q_{ni}(z')] = O(h_n^{d_x}) = o(1)$
        by \eqref{eq:pf-fclt-exp-qni}. Therefore, for any $z, z' \in \scrZ$, $\E[I_{n1}(z,x_0) I_{n1}(z',x_0)] \to R_2(K) \nu(z,z',x_0) p(x_0) \coloneqq \Lambda(z,z',x_0)$ as $n \to \infty$.

    \textbf{Verifying Condition (iii).} Notice that
        \begin{align*}
            &\quad \sum_{i=1}^n \E[Q_{ni}^2]
            = \int_{\RR^{d_x}} K^2(u) \E[H^2(Y_i) \vert X_i=x_0+h_n u] p(x_0+h_n u) d u \nonumber \\
            &\le \int_{\RR^{d_x}} K^2(u) \left(1+\E[H^{2+\gamma}(Y_i) \vert X_i=x_0+h_n u] \right) p(x_0+h_n u) d u \\
            &= \int_{\RR^{d_x}} K^2(u) ( 1+C_F^{2+\gamma}) \left(p(x_0) + h_n u^\intercal \nabla_x p(\bar{x}_0)\right) d u \nonumber \le C_2 R_2(K) + O(h_n) \nonumber < \infty, \nonumber
        \end{align*}
        where
        $C_2 > 0$ is a constant. The first inequality is by the identity that $\E[A^2] = \E[A^2 \ind\{A^2 \le 1\}] + \E[A^2 \ind\{A^2 > 1\}] \le 1 + \E[A^{2+\gamma}]$ for a non-negative random variable $A$ with some constant $\gamma > 0$.

    \textbf{Verifying Condition (iv).} For any $\epsilon > 0$, it holds that
        \begin{align*}
            &\quad \sum_{i=1}^n \E [Q_{ni}^2 \ind\{Q_{ni} > \epsilon\}]
            \le \frac{n\E\left[Q_{ni}^{2+\gamma}\right]}{\epsilon^\gamma}
            = (nh_n^{d_x})^{-\gamma/2} \epsilon^{-\gamma} \E\left[\frac{1}{h_n^{d_x}} K^{2+\gamma}\left(\frac{X_i-x_0}{h_n}\right) H^{2+\gamma}(Y_i)\right] \\
            &= (nh_n^{d_x})^{-\gamma/2} \epsilon^{-\gamma} \int_{\RR^{d_x}} K^{2+\gamma}(u) \E[H^{2+\gamma}(Y_i)|X=x_0+h_n u] p(x_0+h_n u) d u \\
            &\le (nh_n^{d_x})^{-\gamma/2} \epsilon^{-\gamma} C_3 R_2(K) = o(1),
        \end{align*}
        where $C_3 > 0$ is a constant, and the second line is by Markov's inequality.

    \textbf{Verifying Condition (v).} Define $\rho_n^2(z,z') \coloneqq \sum_{i=1}^n \E[(q_{ni}(z)-q_{ni}(z'))^2]$,
        and consider the pseudometric $\rho^2(z,z') \coloneqq R_2(K) \left( \nu(z,z,x_0) - 2\nu(z,z',x_0) + \nu(z',z',x_0) \right) p(x_0) < \infty$.
        It suffices to show that $\rho_n^2(z,z') - \rho^2(z,z') = o(1)$ holds uniformly on $\scrZ$. Indeed, we have
        \begin{align*}
            \rho_n^2(z,z')
            &= \int_{\RR^{d_x}} K^2(u) \left(\nu(z,z,x_0+h_n u) - 2\nu(z,z',x_0+h_n u) + \nu(z',z',x_0+h_n u)\right) p(x_0+h_n u) d u \\
            &= R_2(K) \left(\nu(z,z,x_0) - 2\nu(z,z',x_0) + \nu(z',z',x_0) \right) p(x_0) + O(h_n) \to \rho^2(z,z'),
        \end{align*}
        as $n \to \infty$ for any $z,z' \in \scrZ$.

    Since all conditions hold, $I_{n1}(\cdot,x_0)$ converges to a Gaussian process $\widetilde{\GGG}(\cdot, x_0)$.
    Its finite-dimensional distributions are Gaussian with mean zero and covariance $\Lambda(z,z') = R_2(K) \nu(z,z',x_0) p(x_0)$ for any $z,z' \in \scrZ$. Combining this result with $I_{n2}(\cdot,x_0)$ completes the proof.
    \Halmos \endproof

\subsubsection{Proof of Proposition~\ref{prop:fclt-cont}}
    Notice that
    \begin{align}
        &\quad \sqrt{nh_n^{d_x}} (\hat{f}_n(z,x_0) - f(z,x_0)) = \sqrt{nh_n^{d_x}} \left( \frac{\hat{r}_n(z,x_0)}{\hat{p}_n(x_0)} - \frac{r(z,x_0)}{p(x_0)} \right) \nonumber \\
        &= \sqrt{n h_n^{d_x}} \left(\frac{\hat{r}_n(z,x_0)-r(z,x_0)}{p(x_0)}-\frac{r(z,x_0)}{p^2(x_0)}\left(\hat{p}_{n}(x_0)-p(x_0)\right) + O_{\pr}\left(\frac{1}{nh_n^{d_x}}\right) \right) \label{eq:pf-clt-f-expansion} \\
        &= \frac{\sqrt{n h_n^{d_x}}}{p(x_0)} \left(\hat{r}_n(z,x_0) - f(z,x_0) \hat{p}_{n}(x_0) \right) + o_{\pr}(1) \nonumber \\
        &= \underbrace{\frac{1}{\sqrt{n h_n^{d_x}}} \sum_{i=1}^n \frac{1}{p(x_0)} K\left(\frac{X_i-x_0}{h_n}\right)(F(z;Y_i)-f(z, x_0))}_{\coloneqq J_{n}(z,x_0)} + o_{\pr}(1), \label{eq:pf-clt-f-target}
    \end{align}
    where the third line follows from the first-order Taylor's expansion. In particular, expanding the ratio $A/B$ around $A^*/B^*$ gives $A/B = A^*/B^* + (A-A^*)/B^* - A^*(B-B^*)/(B^*)^2 + O_{\pr}(|A-A^*||B-B^*|+|B-B^*|^2)$. The first two terms in the parentheses of~\eqref{eq:pf-clt-f-expansion} are obtained by setting $A = \hat{r}_n(z,x_0)$, $B = \hat{p}_n(x_0)$, $A^* = r(z,x_0)$ and $B^* = p(x_0)$. By Theorem~2.10 of \cite{PaganUllah99_ec}, $\sqrt{nh_n^{d_x}} (\hat{p}_n(x_0)-p(x_0)) = O_{\pr}(1)$.
    On the other hand, Lemma~\ref{lem:fclt-r} implies that $\sqrt{nh_n^{d_x}} (\hat{r}_n(z,x_0)-r(z,x_0)) = O_{\pr}(1)$ for any $z \in \scrZ$. This justifies the third term in the parentheses of~\eqref{eq:pf-clt-f-expansion}.

    From~\eqref{eq:pf-clt-f-target}, the remaining task is to establish an FCLT for $J_{n}(z,x_0)$.
    Define $s_{ni}(z) \coloneqq \frac{1}{\sqrt{n h_n^{d_x}}} \frac{1}{p(x_0)} K\left(\frac{X_i-x_0}{h_n}\right)(F(z;Y_i)-f(z, x_0))$,
    then $J_{n}(z,x_0) = \sum_{i=1}^n s_{ni}(z)$. By~\eqref{eq:pf-fclt-exp-qni}, the expected value of $s_{ni}(z)$ can be calculated as
    \begin{align*}
        \E[s_{ni}(z)]
        &= \sqrt{\frac{h_n^{d_x}}{n}} \frac{1}{p(x_0)} \left( r(z,x_0)+\frac{h_n^2}{2}\int_{\RR^{d_x}} K(u) u^\intercal \nabla_x^2 r(z,\bar{x}_0) u d u - R_2(K) f(z,x_0) \right).
    \end{align*}
    Define the envelope of $s_{ni}(z)$ as $S_{ni} \coloneqq \frac{1}{\sqrt{n h_n^{d_x}}} \frac{1}{p(x_0)} K\left(\frac{X_i-x_0}{h_n}\right)(H(Y_i)+\E[H(Y_i)|X=x_0])$,
    where $H(\cdot)$ is as defined in the proof of Lemma~\ref{lem:fclt-r}. Proceeding analogously to that proof, we can also verify the five conditions in Theorem 10.6 of \cite{Pollard90_ec} with $q_{ni}(z)$ and $Q_{ni}$ replaced by $s_{ni}(z)$ and $S_{ni}$, respectively. For condition (ii), in particular, it holds that
    \begin{align*}
        \E[J_{n}(z,x_0) J_{n}(z',x_0)] &\to \frac{R_2(K)}{p(x_0)} \E[(F(z;Y) - f(z,x_0)) (F(z';Y) - f(z',x_0))|X=x_0] \\
        &\coloneqq \Psi(z,z',x_0),
    \end{align*}
    for any $z, z' \in \scrZ$. The detailed derivations are omitted for brevity, but are available upon request.

    Hence, $\sqrt{nh_n^{d_x}}(\hat{f}_n(\cdot,x_0) - f(\cdot,x_0))$
    converges in distribution to a Gaussian process $\GGG(\cdot,x_0)$. Its finite-dimensional distributions are uniquely determined, with a mean of zero and a covariance function $\Psi(\cdot,\cdot)$.
\Halmos \endproof

\subsubsection{Completing Proof of Theorem~\ref{thm:clt-cont}}

    Let $C(\scrZ)$ denote the Banach space of continuous functions $\phi(\cdot,x_0): \scrZ \mapsto \mathbb{R}$ endowed with the sup-norm $\|\phi\| = \sup_{z \in \scrZ} |\phi(z,x_0)|$. We work with this space in the ensuing analysis. Note that the FCLT established by Proposition~\ref{prop:fclt-cont} holds in  $C(\scrZ)$. Specifically, the continuity of $f(\cdot,x_0)$ follows from the equicontinuity of $F(z;y)$ in $z$. By definition, for any $z \in \scrZ$ and every $\epsilon > 0$, there exists $\delta > 0$ such that $|F(z;y)-F(z';y)| < \epsilon$ for all $z'$ satisfying $\|z-z'\| < \delta$ and any $y \in \scrY$. When $z$ and $z'$ are sufficiently close so that $\|z-z'\| < \delta$, it follows that $|f(z,x_0)-f(z',x_0)| \le \E_{Y|X=x_0}[|F(z;Y) - F(z';Y)||X=x_0] < \epsilon$, by Jensen's inequality. A similar argument applies to $\hat{f}_n(\cdot,x_0)$ by replacing the probability measure $\pr_{Y|X=x_0}$ with its empirical counterpart $\widehat{\pr}_{Y|X=x_0,n} \coloneqq \sum_{i=1}^n w_n(x_i,x_0) \delta_{y_i}$, where $\delta_{y_i}$ denotes the Dirac point mass at $y_i$. Hence, both $f(\cdot,x_0)$ and $\hat{f}_n(\cdot,x_0)$ are random elements in $C(\scrZ)$.

    Let the min-value function $\vartheta: C(\scrZ) \mapsto \RR$ be defined for all $\phi \in C(\scrZ)$, where $\vartheta(\phi) \coloneqq \inf_{z \in \scrZ} \phi(z,x_0)$. This function is real-valued and measurable with respect to the Borel $\sigma$-algebra induced by the topology of $C(\scrZ)$, which follows from the compactness of $\scrZ$. For any $\phi_1, \phi_2 \in C(\scrZ )$, it holds that $|\vartheta(\phi_1) - \vartheta(\phi_2)|
    = \left|\inf_{z \in \scrZ} \sup_{z' \in \scrZ} \left(\phi_1(z,x_0) -\phi_2(z',x_0)\right)\right|
    \le \sup_{z' \in \scrZ} \left| \phi_1(z',x_0) -\phi_2(z',x_0)\right| = \|\phi_1(z,x_0) - \phi_2(z,x_0)\|$.
    Therefore, the min-value function is 1-Lipschitz continuous.
    By Danskin theorem \citep[Theorem 9.26]{Shapiro21}, $\vartheta(\cdot)$ is directionally differentiable at any point $\varphi(\cdot,x_0) \in C(\scrZ)$, i.e., the limit $\vartheta_{\varphi}'(\varsigma) \coloneqq \lim_{t \downarrow 0} \frac{\vartheta(\varphi + t\varsigma)-\vartheta(\varphi)}{t}$
    exists. It provides a local approximation of $\vartheta(\cdot)$ in the sense that $\vartheta(\varphi + \varsigma) - \vartheta(\varphi) = \vartheta_{\varphi}'(\varsigma) + r(\varsigma)$, where $r(\varsigma)$ denotes the remainder term.
    To ensure that it is negligible, the directional derivative should remain well-defined when $\varsigma$ and $t$ are replaced by sequences $\varsigma_n \to \varsigma$ and $t_n \to t$, respectively. By Proposition~9.72 of \cite{Shapiro21_ec}, $\vartheta(\cdot)$ is also directionally differentiable in the Hadamard sense and  $\vartheta_{\varphi}'(\varsigma) = \inf_{z \in \scrS^*(\varphi)} \varsigma(z,x_0)$,
    where $\scrS^*(\varphi) \coloneqq \argmin_{z \in \scrZ} \varphi (z,x_0)$ is non-empty since $\scrZ$ is compact.

    Observe that $\hat{f}_n(\hat{z}_n(x_0),x_0) = \vartheta(\hat{f}_n(z,x_0))$ and $f^*(x_0) = \vartheta(f(z,x_0))$ with $\scrS^*(f) = \scrZ^*(x_0)$. Invoking the delta theorem \citep[Theorem~9.74]{Shapiro21_ec} with $\varphi = f(\cdot,x_0)$ and $\varsigma = \hat{f}_n(\cdot,x_0) - f(\cdot,x_0)$, we have $\hat{f}_n(\hat{z}_n(x_0),x_0) - f^*(x_0) = \vartheta_{f}'(\hat{f}_n(z,x_0) - f(z,x_0)) + o_{\pr}\left(\frac{1}{\sqrt{nh_n^{d_x}}}\right) = \inf_{z \in \scrZ^*(x_0)} (\hat{f}_n(z,x_0) - f(z,x_0)) + o_{\pr}\left(\frac{1}{\sqrt{nh_n^{d_x}}}\right)$.
    Equivalently, $\sqrt{nh_n^{d_x}}  (\hat{f}_n(\hat{z}_n(x_0),x_0) - f^*(x_0)) = \inf_{z \in \scrZ^*(x_0)} \left(\sqrt{nh_n^{d_x}} (\hat{f}_n(z,x_0) - f(z,x_0))\right) + o_{\pr}(1)$,
    since the directional derivative is positively homogeneous, i.e., $\vartheta_{\varphi}'(t\varsigma) = t\vartheta_{\varphi}'(\varsigma)$ for all $\varsigma \in C(\scrZ)$ and any $t > 0$.

    The proof is completed by applying the FCLT in Proposition~\ref{prop:fclt-cont}, and observing that $\scrZ^*(x_0)$ is a singleton
    and the fact that $\Psi(z,z,x_0) = \msfV(z,x_0)$.
\Halmos \endproof

\subsection{Proof of Theorem~\ref{thm:lin}}
    As a first step, we analyze the error bounds of the budget-constrained wSAA estimator $\hat{f}_n(z^{(m)}(x_0),x_0)$.
    Specifically, the lower bound follows from optimality, whereas the upper bound is obtained by induction using~\eqref{eq:linear-conv}:
    \begin{align} \label{eq:pf-lin-errbound}
        &\quad r(\Gamma) \left( \hat{f}_n(\hat{z}_n(x_0),x_0) - f^*(x_0) \right) \le
        r(\Gamma) \left(\hat{f}_n(z^{(m)}(x_0),x_0) - f^*(x_0)\right) \\
        &\leq \underbrace{r(\Gamma) \left( \hat{f}_n(\hat{z}_n(x_0),x_0) - f^*(x_0) \right)}_{\coloneqq I_1(\Gamma)} + \underbrace{r(\Gamma) \theta^m \left(\hat{f}_n(z^{(0)}_n(x_0),x_0) - \hat{f}_n(\hat{z}_n(x_0),x_0)\right)}_{\coloneqq I_2(\Gamma)},  \nonumber
    \end{align}
    where $r(\Gamma): \RR_+ \mapsto \RR_+$ denotes the scaling factor. We define $T_1 = \frac{r(\Gamma)}{(nh_n^{d_x})^{1/2}}$ and $T_2 = r(\Gamma) \theta^m$.

    \textbf{Part (i).} We set $r(\Gamma) = (\Gamma/\log\Gamma)^{(1-\delta d_x)/2}$. Observe that
    \begin{align*}
        T_1 &= \left(\frac{\kappa^{1-\delta d_x}}{h_0^{d_x}}\right)^{1/2} \left( \frac{\Gamma}{nm} \right)^{(1-\delta d_x)/2} \left( \frac{m- \kappa\log\Gamma}{\kappa\log\Gamma} + 1 \right)^{(1-\delta d_x)/2} \to \left(\frac{\kappa^{1-\delta d_x}}{h_0^{d_x}}\right)^{1/2}, \\
        T_2 &= \left(\frac{1}{\log\Gamma}\right)^{(1-\delta d_x)/2} \Gamma^{(1-\delta d_x)/2 -\kappa\log (1/\theta)} e^{(m - \kappa \log \Gamma)\log \theta} \to 0,
    \end{align*}
    as $\Gamma \to \infty$, since $m \sim \kappa \log\Gamma$ and $(1-\delta d_x)/2 -\kappa\log (1/\theta) \le 0$ when $\kappa \ge (1-\delta d_x)/(2\log (1/\theta))$.
    By Slutsky's theorem and Theorem~\ref{thm:clt-cont}, we have $I_1(\Gamma) \Rightarrow \left(\frac{\kappa^{1-\delta d_x}}{h_0^{d_x}}\right)^{1/2} N(0, \msfV(z^*(x_0),x_0))$
    and $I_2(\Gamma) \xrightarrow{p} 0$
    as $\Gamma \to \infty$.
    Consequently, both bounds in~\eqref{eq:pf-lin-errbound} converge to the same limiting distribution. This completes the proof for part (i).

    \textbf{Part (ii).} We set $r(\Gamma) = \Gamma^{\kappa \log (1/\theta)}$. Similar to part (i), we can show that $T_1 \to 0$ and $T_2 \to 1$
    as $\Gamma \to \infty$, since $\kappa\log (1/\theta) - (1-\delta d_x)/2 < 0$ when  $0< \kappa < (1-\delta d_x)/(2\log (1/\theta))$.
    Then, we have $I_1 \xrightarrow{p} 0$
    and $I_2 \xrightarrow{p} \hat{f}_n(z^{(0)}_n(x_0),x_0) - \hat{f}_n(\hat{z}_n(x_0),x_0) = C_4$
    as $\Gamma \to \infty$, where $C_4>0$ is a constant. This follows from the fact that the continuous function $\hat{f}_n(\cdot,x_0)$ is bounded on the compact set $\scrZ$ under Assumption~\ref{assump:regularity}.
    By a standard $\epsilon$-$\delta$ argument in conjunction
    with the definition of convergence in probability, we can conclude that $r(\Gamma)\left(\hat{f}_n(z^{(m)}_n(x_0),x_0) - f^*(x_0)\right) = O_{\pr}(1)$.
\Halmos \endproof

\subsection{Proof of Theorem~\ref{thm:lin-overopt}}
    Proceeding with the error bounds in~\eqref{eq:pf-lin-errbound}, we set $r(\Gamma) = \Gamma^{(1-\tilde{\kappa})(1-\delta d_x)/2}$. Then,
    \begin{align*}
       \frac{r(\Gamma)}{(nh_n^{d_x})^{1/2}} &= \left(\frac{1}{h_0^{d_x}}\right)^{1/2} \left( \frac{\Gamma}{nm} \right)^{(1-\delta d_x)/2} \left( \frac{m}{\Gamma} \right)^{(1-\delta d_x)/2} \Gamma^{(1-\tilde{\kappa})(1-\delta d_x)/2} \to \left( \frac{c_0^{1-\delta d_x}}{h_0^{d_x}} \right)^{1/2}, \\
        r(\Gamma) \theta^m &= e^{((1-\tilde{\kappa})(1-\delta d_x)/2) \log \Gamma - \log (1/\theta)(m/\Gamma^{\tilde{\kappa}}) \Gamma^{\tilde{\kappa}}} \to 0,
    \end{align*}
    as $\Gamma \to \infty$.
    The remaining steps follow the same reasoning as in the proof of  Theorem~\ref{thm:lin}~(i) and are therefore omitted.
\Halmos \endproof

\subsection{Proof of Theorem~\ref{thm:sublin}}
    We set $r(\Gamma) = \Gamma^{\varsigma(\kappa)}$, where $\varsigma(\kappa) = \min \{(1-\kappa)(1-\delta d_x)/2, \kappa \beta\}$. Similar to Theorem~\ref{thm:lin}, we obtain
    \begin{align} \label{eq:pf-sublin-errbound}
        &\quad r(\Gamma) \left(\hat{f}_n(\hat{z}_n(x_0),x_0) - f^*(x_0)\right) \le r(\Gamma) \left( \hat{f}_n(z^{(m)}_n(x_0),x_0) - f^*(x_0) \right) \\
        &\le \underbrace{r(\Gamma) \left(\hat{f}_n(\hat{z}_n(x_0),x_0) - f^*(x_0)\right)}_{\coloneqq I_1(\Gamma)} + \underbrace{\frac{r(\Gamma)}{m^\beta} \Delta_n(x_0)}_{\coloneqq I_2(\Gamma)}. \nonumber
    \end{align}
    Next, we analyze the limiting behavior of the following two quantities:
    \begin{align*}
        \frac{r(\Gamma)}{(nh_n^{d_x})^{1/2}} &= \left(\frac{1}{h_0^{d_x}}\right)^{1/2} \left(\frac{\Gamma}{nm}\right)^{(1-\delta d_x)/2} \left(\frac{m}{\Gamma^\kappa}\right)^{(1-\delta d_x)/2}  \Gamma^{\varsigma(\Gamma)-(1-\kappa)(1-\delta d_x)/2} \nonumber \\
        &\to \begin{cases}
            \left(c_0^{1-\delta d_x}/h_0^{d_x}\right)^{1/2}, & \text{if~} \kappa \beta \ge (1-\kappa)(1-\delta d_x)/2, \\[-1.2ex]
            0, & \text{o/w},
        \end{cases} \\
        \frac{r(\Gamma)}{m^\beta} &= \left(\frac{\Gamma^\kappa}{m}\right)^\beta \Gamma^{\varsigma(\kappa)-\kappa \beta} \to \begin{cases}
            c_0^{-\beta}, & \text{if~} \kappa \beta \le (1-\kappa)(1-\delta d_x)/2, \\[-1.2ex]
            0, & \text{o/w},
        \end{cases}
    \end{align*}
    as $\Gamma \to \infty$.
    Clearly, the limiting distribution of $\hat{f}_n(z^{(m)}(x_0),x_0)$ is normal if $\frac{r(\Gamma)}{m^\beta}$ asymptotically vanishes. This corresponds to the case when $\kappa \beta \ge (1-\kappa)(1-\delta d_x)/2$. As a result,
    $I_1(\Gamma) \Rightarrow \left(\frac{c_0^{1-\delta d_x}}{h_0^{d_x}}\right)^{1/2} N(0,\msfV(z^*(x_0),x_0))$
    and $I_2(\Gamma) \xrightarrow{p} 0$
    as $\Gamma \to \infty$.
    The conclusion follows from the squeezing theorem, since both bounds in \eqref{eq:pf-sublin-errbound} converge to the same limit.
\Halmos \endproof

\subsection{Proof of Theorem~\ref{thm:superlin}}
    As in the previous proofs, we derive the following error bounds:
    \begin{align}  \label{eq:pf-superlin-errbound}
        &\quad r(\Gamma) \left(\hat{f}_n(\hat{z}_n(x_0),x_0)  -
        f^*(x_0) \right) \le r(\Gamma) \left(\hat{f}_n(z^{(m)}_n(x_0),x_0)  -
        f^*(x_0) \right) \\
        &\le r(\Gamma) \left( \hat{f}_n(\hat{z}_n(x_0),x_0) - f^*(x_0) \right) + \theta^{-1/(\eta-1)} \bigg[ \underbrace{r(\Gamma)^{-\eta^m} \theta^{1/(\eta-1)}  \left( \hat{f}_n(z^{(0)}_n(x_0),x_0) - \hat{f}_n(\hat{z}_n(x_0),x_0) \right)}_{\coloneqq I_2(\Gamma)} \bigg]^{\eta^m} \nonumber \\
        &= \underbrace{r(\Gamma) \left(\hat{f}_n(\hat{z}_n(x_0),x_0)  -
        f^*(x_0) \right)}_{\coloneqq I_1(\Gamma)} + \theta^{-1/(\eta-1)} \bigg[ \underbrace{
        \text{term}~1
        }_{\coloneqq J_1(\Gamma)} + \underbrace{\text{term}~2}_{\coloneqq J_2(\Gamma)} + \underbrace{\text{term}~3}_{\coloneqq J_3(\Gamma)} \bigg]^{\eta^m},
    \end{align}
    where the explicit forms of terms 1 to 3 can be found in Section~\ref{subsec:superlinear}.
    By~\eqref{eq:psi-def}, we can express $J_3(\Gamma) = r(\Gamma)^{\eta^{-m}} \exp(-\psi(z_n^{(0)}(x_0),x_0))$.
    We present an identity to be used in the subsequent analysis:
    \begin{align} \label{eq:pf-superlin-claim}
        \log r(\Gamma)^{\eta^{-m}} &= \exp ((\kappa \log \log \Gamma - m) \log \eta) \exp (-\kappa \log \log \Gamma \log \eta) \log r(\Gamma) \nonumber \\
        &= \Big[ 1 + \big((\kappa \log \log \Gamma - m) \log \eta \big) \exp (\xi)\Big] (\log \Gamma)^{-\kappa \log \eta} \log r(\Gamma) \nonumber \\
        &= (\log \Gamma)^{-\kappa \log \eta} \log r(\Gamma) + o((\log \Gamma)^{-\kappa \log \eta} \log r(\Gamma)),
    \end{align}
    where the second line is by the mean value theorem.

    \textbf{Part (i).} We set $r(\Gamma) = (\Gamma/\log\log \Gamma)^{(1-\delta d_x)/2}$. It follows that
    \begin{align*}
        \frac{r(\Gamma)}{(nh_n^{d_x})^{1/2}} = \left(\frac{\kappa^{1-\delta d_x}}{h_0^{d_x}}\right)^{1/2} \left(\frac{\Gamma}{nm} \left(\frac{m - \kappa\log\log\Gamma}{\kappa\log\log\Gamma} + 1 \right)\right)^{(1-\delta d_x)/2} \to \left(\frac{\kappa^{1-\delta d_x}}{h_0^{d_x}}\right)^{1/2},
    \end{align*}
    as $\Gamma \to \infty$.
    Therefore, $I_1(\Gamma) \Rightarrow \left(\frac{\kappa^{1-\delta d_x}}{h_0^{d_x}}\right)^{1/2} N(0,\msfV(z^*(x_0),x_0))$
    as $\Gamma \to \infty$.
    Next, we address $J_1(\Gamma)$, $J_2(\Gamma)$ and $J_3(\Gamma)$ appearing in the upper bound. Substituting $r(\Gamma)$ into~\eqref{eq:pf-superlin-claim} yields
    \begin{align} \label{eq:pf-superlin-logr-i}
        \log r(\Gamma)^{\eta^{-m}}
        &= \frac{1}{2} (1-\delta d_x) (\log \Gamma)^{1-\kappa \log \eta} \left( 1 - \frac{\log\log\log\Gamma}{\log\Gamma} \right) + o((\log \Gamma)^{1-\kappa \log \eta}),
    \end{align}
    which is bounded since $1-\kappa \log \eta \le 0$. Therefore, $r(\Gamma)^{\eta^{-m}} = O(1)$, so that $\frac{r(\Gamma)^{\eta^{-m}}}{(nh_n^{d_x})^{1/2}} \to 0$
    as $\Gamma \to \infty$. Consequently, $J_1(\Gamma)\xrightarrow{p} 0$ and $J_2(\Gamma) \xrightarrow{p} 0$
    as $\Gamma \to \infty$.
    It remains to show that $J_3(\Gamma)^{\eta^m} \xrightarrow{p} 0$ as $\Gamma \to \infty$. Equivalently, we need to verify that
    \begin{align} \label{eq:pf-superlin-J3-conv}
        \lim_{\Gamma \to \infty} \pr \left(\frac{J_3(\Gamma)}{\epsilon^{\eta^{-m}}} > 1 \right) + \pr \left(\frac{J_3(\Gamma)}{\epsilon^{\eta^{-m}}} < -1 \right) = 0, \quad \forall \epsilon > 0.
    \end{align}
    Clearly, the denominator $\epsilon^{\eta^{-m}} \to 1$
    as $\Gamma \to \infty$, since $\log \epsilon^{\eta^{-m}} = \eta^{-m} \log \epsilon \to 0$
    as $m \to \infty$, which holds for any $\eta > 1$ and $\epsilon > 0$. For the numerator, it follows from \eqref{eq:pf-superlin-logr-i} that
    \begin{align}
        \log J_3(\Gamma) = -\psi(z^{(0)}_n(x_0),x_0) + \frac{1}{2} (1-\delta d_x) (\log \Gamma)^{1-\kappa \log \eta} \left( 1 - \frac{\log\log\log\Gamma}{\log\Gamma} \right)+ o(1). \label{eq:pf-superlin-logJ3}
    \end{align}
    \textbf{Case (a).} If $\kappa > 1/\log \eta$, then~\eqref{eq:pf-superlin-logJ3} simplifies to $\log J_3(\Gamma) = -\psi(z^{(0)}_n(x_0),x_0) + o(1)$.
    Since $\psi(z^{(0)}_n(x_0),x_0) > 0$, then $0 < \exp(-\psi(z^{(0)}_n(x_0),x_0)) < 1$.
    Combining this with the fact that $\epsilon^{\eta^{-m}} \to 1$ as $\Gamma \to \infty$, we obtain
    \begin{align*}
        \frac{J_3(\Gamma)}{\epsilon^{\eta^{-m}}} = \frac{\exp(-\psi(z^{(0)}_n(x_0),x_0) + o(1))}{1+o(1)} \rightarrow \exp(-\psi(z^{(0)}_n(x_0),x_0)) \in (0,1),
    \end{align*}
    as $\Gamma \to \infty$. Hence, \eqref{eq:pf-superlin-J3-conv} holds.
    \textbf{Case (b).} If $\kappa = 1/\log \eta$, then \eqref{eq:pf-superlin-logJ3} becomes
    \begin{align}
        \log J_3(\Gamma) = - \psi(z^{(0)}_n(x_0),x_0) + \frac{1}{2} (1-\delta d_x) \left( 1 - \frac{\log\log\log\Gamma}{\log\Gamma} \right) + o(1). \label{eq:pf-superlin-logJ3-b}
    \end{align}
    \textbf{Sub-case (1).} When $\psi(z^{(0)}_n(x_0),x_0) > (1-\delta d_x)/2$, $J_3(\Gamma) \xrightarrow{p} \exp\left(- \psi(z^{(0)}_n(x_0),x_0) + \frac{1}{2} (1-\delta d_x) \right)$ $\in (0,1)$
    as $\Gamma \to \infty$. \textbf{Sub-case (2).} When $\psi(z^{(0)}_n(x_0),x_0) = (1-\delta d_x)/2$, \eqref{eq:pf-superlin-logJ3-b} reduces to $\log J_3(\Gamma) = - \frac{1}{2}(1-\delta d_x) \frac{\log \log \log \Gamma}{\log \Gamma} + o(1) \to 0$
    as $\Gamma \to \infty$. Since $\log \epsilon^{\eta^{-m}} \to 0$ as $\Gamma \to \infty$,
    a more refined analysis
    is warranted. According to~\eqref{eq:pf-superlin-claim}, we have
    \begin{align} \label{eq:pf-superlin-logeps-finer}
        \log \epsilon^{\eta^{-m}} &= \exp((\kappa \log \log \Gamma -m) \log \eta) \exp(-\kappa \log \log \Gamma \log \eta) \log \epsilon = (1 + o(1)) (\log \Gamma)^{-\kappa \log \eta} \log \epsilon \nonumber \\
        &= (\log \Gamma)^{-\kappa \log \eta} \log \epsilon  + o( (\log \Gamma)^{-\kappa \log \eta}) = (\log \Gamma)^{-1} \log \epsilon  + o(1).
    \end{align}
    This indicates that $\epsilon^{\eta^{-m}}$ is dominated by $J_3(\Gamma)$. As a result, the quantity $\log \left( \frac{J_3(\Gamma)}{\epsilon^{\eta^{-m}}} \right) =  - \frac{1}{2}(1-\delta d_x) \frac{\log \log \log \Gamma}{\log \Gamma} + O((\log \Gamma)^{-1})$
    is negative with sufficiently large $\Gamma$,
    confirming that \eqref{eq:pf-superlin-J3-conv} also holds for this boundary case.
    Combining the two sub-cases regarding $\psi(z^{(0)}_n(x_0),x_0)$, we complete the proof of part (i).

    \textbf{Part (ii).} We set $r(\Gamma) = \exp(\psi(z^{(0)}_n(x_0),x_0) (\log \Gamma)^{\kappa \log \eta})$.
    It follows that
    \begin{align} \label{eq:pf-superlin-ii}
        0 < \frac{r(\Gamma)^{\eta^{-m}}}{(nh_n^{d_x})^{1/2}}
        \le \frac{\exp(\psi(z^{(0)}_n(x_0),x_0) (\log \Gamma)^{\kappa \log \eta})}{(nh_n^{d_x})^{1/2}} = \frac{r(\Gamma)}{(nh_n^{d_x})^{1/2}},
    \end{align}
    where the second inequality holds for sufficiently large $\Gamma$, since $\eta^{-m} < 1$ whenever $\eta > 1$. In what follows, we focus on $\frac{r(\Gamma)}{(nh_n^{d_x})^{1/2}}$ and show that it goes to 0 as $\Gamma \to \infty$.
    Notice that
    \begin{align} \label{eq:pf-superlin-ratio}
        & \frac{r(\Gamma)}{(nh_n^{d_x})^{1/2}} =
        \exp(\psi(z^{(0)}_n(x_0),x_0) (\log \Gamma)^{\kappa \log \eta} - \frac{1}{2}(1-\delta d_x) \log \Gamma) \nonumber \\
        &\qquad \cdot \left(\frac{\Gamma}{nm}\right)^{(1-\delta d_x)/2} \left(\frac{m}{\kappa \log \log \Gamma}\right)^{(1-\delta d_x)/2} \left(\frac{\kappa^{1-\delta d_x}}{h_0^{d_x}}\right)^{1/2} (\log \log \Gamma)^{(1-\delta d_x)/2} \to 0,
    \end{align}
    as $\Gamma \to \infty$, provided that the condition $\psi(z^{(0)}_n(x_0),x_0) (\log \Gamma)^{\kappa \log \eta} - \frac{1}{2}(1-\delta d_x) \log \Gamma < 0$ is satisfied for sufficiently large $\Gamma$.
    This is indeed true if either of the following cases applies:
    \textbf{Case (a).} If $\kappa < 1/\log \eta$, then $(\log \Gamma)^{\kappa \log \eta}$ grows more slowly than $\log\Gamma$, so the negative part dominates.
    \textbf{Case (b).} If $\kappa = 1/\log \eta$ and $\psi(z^{(0)}_n(x_0),x_0) < (1-\delta d_x)/2$, then the required condition is clearly satisfied.
    With~\eqref{eq:pf-superlin-ratio} established, it follows that $I_1(\Gamma) \xrightarrow{p} 0$
    as $\Gamma \to \infty$.
    Reviewing ~\eqref{eq:pf-superlin-ii}, we can further deduce that $\frac{r(\Gamma)^{\eta^{-m}}}{(nh_n^{d_x})^{1/2}} \to 0$
    as $\Gamma \to \infty$.
    As a consequence,
    $J_1(\Gamma)\xrightarrow{p} 0$ and $J_2(\Gamma) \xrightarrow{p} 0$
    as $\Gamma \to \infty$.
    Finally, we inspect the term $J_3(\Gamma)$. Referring to~\eqref{eq:pf-superlin-claim}, we have
    \begin{align*}
        \log J_3(\Gamma)
        &= -\psi(z^{(0)}_n(x_0),x_0) + (\log \Gamma)^{-\kappa \log \eta} \left(\psi(z^{(0)}_n(x_0),x_0) (\log \Gamma)^{\kappa \log \eta} \right) + o(1) \to 0
    \end{align*}
    as $\Gamma \to \infty$.
    For any given $\epsilon > 0$, $\log \left( \frac{J_3(\Gamma)}{\epsilon^{\eta^{-m}}}\right) = -\log \epsilon  (\log \Gamma)^{-\kappa \log \eta} +  o(1)$
    by \eqref{eq:pf-superlin-logeps-finer}. Therefore, there exists $\epsilon > 1$ such that $\log \left( \frac{J_3(\Gamma)}{\epsilon^{\eta^{-m}}}\right) < 0$ when $\Gamma$ is large enough, i.e., $J_3(\Gamma)/\epsilon^{\eta^{-m}} \in (0,1)$ for all sufficiently large $\Gamma$,
    which in turn implies that $J_3(\Gamma)^{\eta^m} = O(1)$.
    To summarize, the lower bound in \eqref{eq:pf-superlin-errbound} is $o_{\pr}(1)$, whereas the corresponding upper bound is $O_{\pr}(1)$. Hence, we conclude that $r(\Gamma) ( \hat{f}_n(z^{(m)}_n(x_0),x_0) - f^*(x_0) ) = O_{\pr}(1)$.
\Halmos \endproof

\subsection{Proof of Theorem~\ref{thm:superlin-overopt}}
    We set $r(\Gamma) = (\Gamma/\log \Gamma)^{(1-\delta d_x)/2}$. Following similar reasoning as in~\eqref{eq:pf-superlin-claim}, we obtain
    \begin{align*}
        \log r(\Gamma)^{\eta^{-m}}
        &= (1 + o(1))  \Gamma^{-\tilde{\kappa} \log \eta} \frac{1}{2}(1-\delta d_x) (\log \Gamma - \log \log \Gamma) = O(\Gamma^{-\tilde{\kappa} \log \eta} \log \Gamma) = o(1),
    \end{align*}
    where the last equality follows immediately from $\tilde{\kappa} > 0$ and $\eta > 1$. Therefore, $r(\Gamma)^{\eta^{-m}} \to 1$ as $\Gamma \to \infty$, which implies that $\frac{r(\Gamma)^{\eta^{-m}}}{(nh_n^{d_x})^{1/2}} \to 0$
    as $\Gamma \to \infty$, and thus
    $J_1(\Gamma)\xrightarrow{p} 0$ and $J_2(\Gamma) \xrightarrow{p} 0$
    as $\Gamma \to \infty$. Moreover, we have $\log J_3(\Gamma) = -\psi(z^{(0)}_n(x_0),x_0) + \log r(\Gamma)^{\eta^{-m}} \to -\psi(z^{(0)}_n(x_0),x_0) < 0$
    as $\Gamma \to \infty$. Therefore, $J_3(\Gamma)^{\eta^m} \to C_5$ as $n \to \infty$, where $C_5 \in (0,1)$ is a constant.
    For every $\epsilon > 0$, $\lim_{\Gamma \to \infty} \pr( |I_2(\Gamma)^{{\eta^m}}| > \epsilon ) \le \lim_{\Gamma \to \infty} \big( \pr( |J_1(\Gamma)| > 1) + \pr( |J_2(\Gamma)| > 1) + \pr( |J_3(\Gamma)| > 1) \big) = 0$
    by the identity that $\epsilon^{\eta^{-m}} \to 1$ as $\Gamma \to \infty$. Hence, we have $I_2(\Gamma)^{{\eta^m}} \xrightarrow{p} 0$ as $\Gamma \to \infty$.
    On the other hand,
    \begin{align*}
        \frac{r(\Gamma)}{(nh_n^{d_x})^{1/2}} = \left(\frac{\tilde{\kappa}^{1-\delta d_x}}{h_0^{d_x}}\right)^{1/2} \left(\frac{\Gamma}{nm} \left(\frac{m - \tilde{\kappa}\log\Gamma}{\tilde{\kappa}\log\Gamma} + 1 \right)\right)^{(1-\delta d_x)/2} \to \left(\frac{\tilde{\kappa}^{1-\delta d_x}}{h_0^{d_x}}\right)^{1/2},
    \end{align*}
    as $\Gamma \to \infty$, which leads to $I_1(\Gamma) \Rightarrow \left(\frac{\tilde{\kappa}^{1-\delta d_x}}{h_0^{d_x}}\right)^{1/2} N(0,\msfV(z^*(x_0),x_0))$
    as $\Gamma \to \infty$. With both bounds in~\eqref{eq:pf-superlin-errbound} converging to the same limiting distribution, the proof is complete.
\Halmos \endproof

\subsection{Proof of Corollary~\ref{cor:ci-cvg}}

\subsubsection{A Technical Lemma}
\begin{proposition} \label{prop:cond-var-consistency}

Suppose Assumptions~\ref{assump:regularity}--\ref{assumption:band} hold, and
    $\E[|F^2(z;Y)| (\log |F^2(z;Y)|)^{+}] < \infty$ for every $z \in \scrZ$.
    Then, $\hat{\sigma}_n^2(z,x_0) \xrightarrow{a.s.} \sigma^2(z,x_0)$
    as $n\to\infty$, for all $x_0 \in \scrX$ and every $z \in \scrZ$,
    where $\hat{\sigma}_n^2(z,x_0)$ is defined in \eqref{eq:sigma_hat}.
\end{proposition}

\proof{Proof of Proposition~\ref{prop:cond-var-consistency}.}
    We first fix a $z \in \scrZ$. Clearly, $\E[|F(z;Y)| (\log |F(z;Y)|)^{+}] < \infty$. Applying Theorem~3 of \cite{Walk10_ec} to the random variable $Y_z \coloneqq F(z;Y)$, we obtain $\hat{f}_n(z,x_0) \xrightarrow{a.s.} f(z,x_0)$ as $n\to\infty$, for all $x_0 \in \scrX$. By Lemma EC.7 of \cite{BertsimasKallus20_ec}, this result can be generalized to every $z \in \scrZ$. Specifically, $\hat{f}_n(z,x_0) \xrightarrow{a.s.} f(z,x_0)$ as $n\to\infty$ for all $x_0 \in \scrX$ and every $z \in \scrZ$. The pointwise convergence of $\hat{f}_n(z,x_0)$ further implies uniform convergence \citep[Lemma EC.5]{BertsimasKallus20_ec}, i.e., $\sup_{z \in \scrZ} \left|\hat{f}_n(z,x_0)-f(z,x_0)\right| \xrightarrow{a.s.} 0$
    as $n\to\infty$.

    We write $\hat{\sigma}_n^2(z,x_0) = \sum_{i=1}^n w_n(x_i,x_0)F^2(z;y_i) - \left( \sum_{i=1}^n w_n(x_i,x_0)F(z;y_i) \right)^2 \coloneqq \hat{s}_{n}^2(z,x_0) - \hat{f}_n^2(z,x_0)$.
    Since $\hat{f}_n(\cdot,x_0)$ is a uniformly convergent sequence of bounded functions on the compact set $\scrZ$, there exists a constant $C_6 > 0$ such that $\sup_{z \in \scrZ}|\hat{f}_n(z,x_0)| \le C_6$. The limiting function $f(\cdot,x_0)$ is also uniformly bounded by this constant.
    Note that $C_6$ is essentially $C_f$ under Assumption~\ref{assump:smooth-1}. It then follows that $\sup_{z \in \scrZ} \left|\hat{f}_n^2(z,x_0) - f^2(z,x_0)\right| = \sup_{z \in \scrZ} \left|\hat{f}_n(z,x_0) - f(z,x_0)\right| \left|\hat{f}_n(z,x_0) + f(z,x_0)\right| \le 2C_6\sup_{z \in \scrZ} \left|\hat{f}_n(z,x_0) - f(z,x_0)\right| \xrightarrow{a.s.} 0$
    as $n\to\infty$.
    Hence, $\hat{f}_n^2(z,x_0)$ converges to $f^2(z,x_0)$ uniformly over $z \in \scrZ$.

    The remaining task is to show that, for all $x_0 \in \scrX$, $\hat{s}_{n}^2(z,x_0) \xrightarrow{a.s.} s^2(z,x_0) \coloneqq \E[F^2(z;Y)|X=x_0]$ as $n\to\infty$, uniformly over $z \in \scrZ$. Again, we first fix a $z \in \scrZ$ and define a random variable $Y_z' \coloneqq F^2(z;Y)$. By assumption, $\E[|Y_z'| (\log |Y_z'|)^{+}] < \infty$ for all $z \in \scrZ$. Applying Theorem 3 of \cite{Walk10_ec} to $Y_z'$ gives $\hat{s}_{n}^2(z,x_0) \xrightarrow{a.s.} s^2(z,x_0)$ as $n\to\infty$, for all $x_0 \in \scrX$.
    In the next step, we extend this pointwise convergence to hold uniformly in $z$. Consider the set $\scrZ' \coloneqq \scrZ \cap \QQ^{d_z} \cup S$, where $S \coloneqq \{z \in \scrZ: \scrV_\epsilon(z) \cap \scrZ = \{z\}, \forall \epsilon > 0\}$ denotes the set of isolated points in $\scrZ$. For each $z \in S$, its neighborhood $\scrV_\epsilon(z) \coloneqq \{z' \in \mathrm{int}(\scrZ):\|z'-z\|_{\infty} \leq \epsilon, \epsilon>0 \}$ contains no points of $\scrZ$ except $z$ itself, where $\mathrm{int}(\scrZ)$ denotes the interior of $\scrZ$. Therefore, $\scrZ'$ is countable and dense by construction. For each fixed $z' \in \scrZ'$, define $A(z') \coloneqq \left\{\omega \in \Omega: \lim_{n \to \infty} \hat{s}_{n}^2(z',x_0)(\omega) = s^2(z',x_0)\right\}$, where $\omega$ represents a realization of the dataset $\scrD_n$.
    By the definition of almost‑sure convergence, each set $A(z')$ has probability one. Since the probability measure is continuous and $\scrZ$ is countable, we have $\pr \left(\bigcap_{z' \in \scrZ'} A(z')\right) = 1$.
    We take a particular $\omega$ for which the event $A(z')$ occurs. By the triangle inequality, we have $\left|\hat{s}_n^2(z,x_0) - s^2(z,x_0)\right| \le \underbrace{\left|\hat{s}_n^2(z,x_0) - \hat{s}_n^2(z',x_0)\right|}_{\coloneqq V_1} + \underbrace{\left|s^2(z,x_0) - s^2(z',x_0)\right|}_{\coloneqq V_2} + \underbrace{\left|\hat{s}_n^2(z',x_0) - s^2(z',x_0)\right|}_{\coloneqq V_3}$.
    Let $\epsilon > 0$ be given.
    Since $F(z;y)$ is equicontinuous in $z$, there exists $\delta > 0$ such that $|F(z;y)-F(z';y)| \le \epsilon/3$ whenever $\|z-z'\| \le \delta$. It then follows that $V_1 \le \sum_{i=1}^n w_{n}(x_i,x_0)|F^2(z;y_i)-F^2(z';y_i)| \le 2C_7 |F(z;y)-F(z';y)| = \epsilon'/3$,
    where $\epsilon' > \epsilon$ and $C_7 = \sup_{z \in \scrZ, y \in \scrY}|F(z;y)|$ is some constant under Assumption~\ref{assump:envelope}. Similarly, we can show that $V_2 \le \E[|F^2(z;Y)-F^2(z';Y)||X=x_0] \le \epsilon'/3$.
    Since $\hat{s}_{n}^2(z,x_0) \xrightarrow{a.s.} s^2(z,x_0)$ as $n\to\infty$ for every fixed $z \in \scrZ$, there exists $k \in \N_+$ such that $V_3 \le \epsilon'/3$,
    for all $n \ge k$. Collecting the bounds for $V_1, V_2$ and $V_3$, we conclude that $|\hat{s}_n^2(z,x_0) - s^2(z,x_0)| \le \epsilon'$.
    Here, $\epsilon$ can be made arbitrarily small, and so can $\epsilon'$. Hence, the above argument holds for all $z \in \scrZ$. As the set of $\omega$ for which this result holds has probability one, the proof is complete.
\Halmos \endproof

\subsubsection{Completing Proof of Corollary~\ref{cor:ci-cvg}}
    We first show that for all $x_0 \in \scrX$, $\hat{\sigma}_n^2(\hat{z}_n(x_0),x_0)$ is a strongly consistent estimator of $\sigma^2(z^*(x_0),x_0)$, namely $\sigma_n^2(\hat{z}_n(x_0),x_0) \xrightarrow{a.s.} \sigma^2(z^*(x_0),x_0)$ as $n\to\infty$ for every $\hat{z}_n(x_0) \in \hat{\scrZ}(x_0)$, where $\hat{\scrZ}_n(x_0)$ denotes the set of optimal solutions to the wSAA problem~\eqref{eq:cso-wSAA}. Let $\epsilon > 0$ be given. By Proposition~\ref{prop:cond-var-consistency}, there exists $k_1 \in \N_+$ such that $W_1 \coloneqq \left|\hat{\sigma}_n^2(\hat{z}_n(x_0),x_0) - \sigma^2(\hat{z}_n(x_0),x_0)\right| \le \epsilon/2$
    for all $n \ge k_1$. Since $\scrZ^*(x_0) = \{z^*(x_0)\}$ is a singleton, then $\hat{z}_n(x_0) \xrightarrow{a.s.} z^*(x_0)$ as $n\to\infty$ for every $\hat{z}_n(x_0) \in \hat{\scrZ}_n(x_0)$ \citep[Theorem~6]{BertsimasKallus20_ec}. It is straightforward to verify that $F^2(z;y)$ is also equicontinuous in $z$. As a result, $s^2(\cdot,x_0)$ is continuous for all $x_0 \in \scrX$. Since $f^2(\cdot,x_0)$ is continuous for all $x_0 \in \scrX$ as well, the continuity of $\sigma^2(\cdot,x_0) = s^2(\cdot,x_0)-f^2(\cdot,x_0)$ follows immediately for all $x_0 \in \scrX$. Therefore, there exists $k_2 \in \N_+$ such that $W_2 \coloneqq \left|\sigma^2(\hat{z}_n(x_0),x_0) - \sigma^2(z^*(x_0),x_0)\right| \le \epsilon/2$
    for all $n \ge k_2$. For all $x_0 \in \scrX$, it then follows that $\left|\hat{\sigma}_n^2(\hat{z}_n(x_0),x_0) - \sigma^2(z^*(x_0),x_0)\right| \le W_1 + W_2 \le \epsilon/2 + \epsilon/2 \le \epsilon$,
    for every $\hat{z}_n(x_0) \in \scrZ(x_0)$ and all $n \ge \max(k_2,k_3)$. Sending $\epsilon$ to zero, we obtain the desired result.

    It remains to show that for all $x_0 \in \scrX$, $nh_n^{d_x} \sum_{i=1}^n w_n^2(x_i,x_0) \xrightarrow{p} \frac{R_2(K)}{p(x_0)}$ as $n \to \infty$.
    Observe that $nh_n^{d_x} \sum_{i=1}^n w_n^2(x_i,x_0) = \frac{(nh_n^{d_x})^{-1} \sum_{i=1}^n K^2((x_i-x_0)/h_n)}{\left((nh_n^{d_x})^{-1} \sum_{i=1}^n K((x_i-x_0)/h_n)\right)^2} \coloneqq \frac{T_1}{T_2}$,
    where $T_2 = p^2(x_0)+o_{\pr}(1)$ by Theorem~2.6 of \cite{PaganUllah99_ec}. From Bochner's lemma \citep[Appendix~A.2.6]{PaganUllah99_ec}, for any $s \ge 0$, $\frac{1}{h_n^{d_x}} \E\left[K^s\left(\frac{X-x_0}{h_n}\right)\right] = p(x_0)R_s(K) + o(1)$,
    where $R_s(K) \coloneqq \int_{\RR^{d_x}} K^s(u) du$ is the $s$-th power integral of the kernel function. It then follows that $\E[T_1] = p(x_0)R_2(K) + o(1)$
    and $\vari(T_1) = \frac{1}{nh_n^{d_x}} \left(\frac{1}{h_n^{d_x}}\E\left[ K^4\left(\frac{X-x_0}{h_n}\right)\right]\right) - \frac{1}{n} \left(\frac{1}{h_n^{d_x}}\E\left[ K^2\left(\frac{X-x_0}{h_n}\right)\right] \right)^2 = \frac{1}{nh_n^{d_x}} \left(p(x_0) R_4(K) + o(1) \right) - \frac{1}{n} \left(p(x_0) R_2(K) + o(1) \right)^2 = o(1)$,
    where $R_4(K) < \infty$ by assumption.
    For any given $\epsilon>0$, $\pr(|T_1-p(x_0)R_2(K) | \ge \epsilon) \le \pr(|T_1 - \E[T_1]| \ge \epsilon/2) + \pr(|\E[T_1] - p(x_0)R_2(K)| \ge \epsilon/2) \le 4\vari(T_1)/\epsilon^2 + \pr(|\E[T_1] - p(x_0)R_2(K)| \ge \epsilon/2) \to 0$
    as $n \to \infty$, by the Chebyshev's inequality. Hence, we conclude that $\frac{T_1}{T_2} = \frac{R_2(K)}{p(x_0)} + o_{\pr}(1)$.

    Finally, $\frac{\sqrt{nh_n^{d_x}}(\hat{f}_n(\hat{z}_n(x_0),x_0) - f^*(x_0))}{\sqrt{\sigma^2(z^*(x_0),x_0)R_2(K)/p(x_0)}} = \frac{\hat{f}_n(\hat{z}_n(x_0),x_0) - f^*(x_0)}{\sqrt{\hat{\sigma}_n^2(\hat{z}_n(x_0),x_0) \sum_{i=1}^n w_n^2(x_i,x_0)}} \sqrt{\frac{\hat{\sigma}_n^2(\hat{z}_n(x_0),x_0)}{\sigma^2(z^*(x_0),x_0)}} \sqrt{\frac{nh_n^{d_x} \sum_{i=1}^n w_n^2(x_i,x_0)}{R_2(K)/p(x_0)}} \Rightarrow N(0,1)$ as $n \to \infty$.
\Halmos \endproof

\subsection{Proof of Corollary~\ref{cor:contrained-CI}}
    According to our assumption on the algorithm $\scrA$, we have $\hat{f}_n(z^{(m)}_n(x_0),x_0) \to \hat{f}_n(\hat{z}_n(x_0),x_0)$ as $m \to \infty$. Recall Assumption~\ref{assump:uni-sol}, $z^*(x_0)$ is unique.
    Then given the continuity of $\hat{f}_n(\cdot,x_0)$, we can show, by contradiction, that $z^{(m)}_n(x_0) \to \hat{z}_n(x_0)$ as $m \to \infty$. The proof is then completed by joining Corollary~\ref{cor:ci-cvg} with Theorems~\ref{thm:lin}-\ref{thm:superlin-overopt}.
 \Halmos \endproof

\section{Newsvendor Problem with Sublinearly Convergent Algorithm} \label{app:exp-nv}

\begin{figure}[ht]
    \FIGURE{\includegraphics[width=\textwidth]{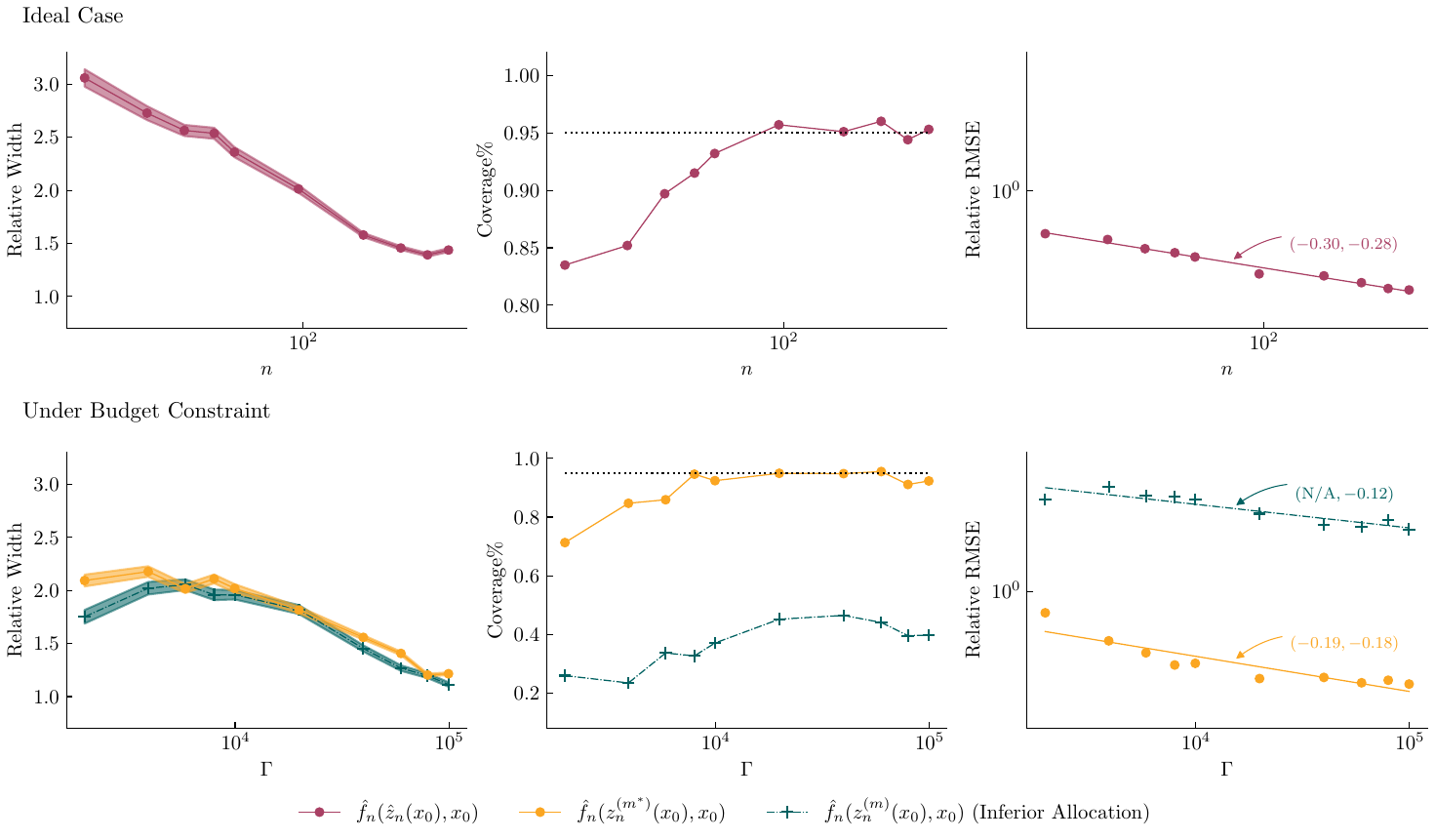}}
{Newsvendor Problem. \label{fig:exp-NV}}
{The first number in each pair of annotated parentheses represents the theoretical convergence rate of relative RMSE, which is $n^{-(1-\delta d_x)/2}$ for $\hat{f}_n(\hat{z}_n(x_0),x_0)$ (Theorem~\ref{thm:clt-cont}) and $\Gamma^{-\kappa^*\beta}$ for $\hat{f}_n(z^{(m^*)}_{n}(x_0),x_0)$ (Theorem~\ref{thm:sublin}). The second number indicates the empirical slope obtained from regressing log relative RMSEs on $\log n$ or $\log \Gamma$.  ``N/A'' denotes cases where the theoretical rate is not applicable. The optimization algorithm is sublinearly convergent.
}
\end{figure}

Consider a single-item newsvendor problem where we determine ordering quantities for a perishable good under uncertain demand.
Each unit of lost sales incurs an underage cost of $c_u=10$, while each unit of excess inventory incurs an overage cost of $c_o=2$.
The total inventory cost for decision $z$ and random demand $Y$ is $F(z; Y) = c_u (Y-z)^+ + c_o(z-Y)^+$.
We assume the covariate $X=(X^1,X^2)$ is two-dimensional with independent components where $X^1 \sim N(20,2^2)$ and $X^2 \sim \mathsf{LogNorm}(1,0.3^2)$.
The new covariate observation $x_0=(x_0^1, x_0^2)$ is set as the 25\% quantile ($\tau=0.25$) of the marginal distributions of $X^1$ and $X^2$, respectively.
Furthermore, we assume
the conditional distribution of $Y$ given $X$ is a left-truncated normal distribution at zero where the pre-truncation normal has a mean of $100+(X^1 - 20) + X^2 ( 2\ind_{(-\infty,2]} + 4 \ind_{(2,4]} + 6\ind_{(4,6]} + 8\ind_{(6,\infty)})$ and a standard deviation of 3.
Since $\hat{f}_n(z,x_0)$ is a weighted sum of convex functions $F(\cdot;y_i)$ for $i \in [n]$, it is convex and  $L$-Lipschitz continuous with $L = \max(c_u, c_o)$.
To solve the resulting wSAA problem, we use projected subgradient descent method,
which converges sublinearly with parameter $\beta = 1/2$ as defined in \eqref{eq:sublinear-conv}.
The subgradient of $\hat{f}_n(z,x_0)$ has the explicit form
$\partial_z \hat{f}_n(z,x_0) = \sum_{i=1}^n w_n(x_i,x_0) \partial_z F(z;y_i)$, where
where $\partial_z F(z;y_i) = \{-c_u\}$ for $z < y_i$, $\partial_z F(z;y_i) = \{-c_o\}$ for $z > y_i$, and $\partial_z F(z;y_i) =  [-c_u, c_o]$ for $z=y_i$.
Following the bike-sharing example in Section~\ref{sec:exper},
we use the Gaussian kernel function $K(u) = \exp(-\|u\|^2/2)$, with bandwidth exponent $\delta = 1/(d_x+3)=1/5$. We evaluate the performance of the confidence intervals~\eqref{eq:ci-wSAA} and \eqref{eq:ci-wSAA-constrained}.
The results are presented in Figure~\ref{fig:exp-NV}.

The convergence rate of the relative RMSE of the wSAA estimator $\hat{f}_n(\hat{z}_n(x_0),x_0)$ without computational budget constraints aligns with the theoretical value in Theorem~\ref{thm:clt-cont}. As the sample size $n$ increases, the widths of confidence intervals shrink while their coverages approach the target level when $n$ is larger than $10^2$. Similar patterns emerge for the budget-constrained wSAA estimator $\hat{f}_n(z^{(m^*)}_{n}(x_0),x_0)$ under the optimal budget allocation specified in Theorem~\ref{thm:sublin}.
In addition, we compare the optimal budget allocation to an inferior allocation that fixes algorithm iterations at 50 regardless of computational budget.
This inferior allocation fails to achieve the desired relative RMSE convergence rate and  results in poor confidence interval coverage. Such a heuristic allocation with fixed algorithm iterations  cannot effectively manage the statistical-computational tradeoff.

\section{High-order Polynomial Function with Superlinearly Convergent Algorithm} \label{app:exp-fourthpoly}

Consider a fourth-order polynomial cost function $F(z;Y) = \sum_{j=1}^{d_z} a^j(z^j-b^j Y^j)^4$,
where $Y\in \RR^{d_y}$  and $z \in \RR^{d_z}$ with $d_y = d_z = 2$.
The coefficients $a^j$'s and $b^j$'s are generated from $N(20,15)$ and $\mathsf{Unif}[-5,-1]$, respectively.
We assume the covariate $X=(X^1,X^2)$ is two-dimensional with independent components, where $X^1 \sim N(10,4)$ and $X^2 \sim N(8,1)$.
The conditional distribution of $Y$ given $X$ is multivariate normal with mean vector $(\log(X^1+4)+5, \sqrt{|X^2|}+10)^\intercal$ and covariance being an identity matrix.
To solve the wSAA problem,
we use Newton's method with Armijo backtracking parameters $a = 0.1$ and $b = 0.9$.
The gradient and Hessian of $\hat{f}_n(z,x_0)$ are computed as $\nabla_z \hat{f}_n(z,x_0) = \sum_{i=1}^n w_n(x_i,x_0) \nabla_z F(z;y_i)$ and $\nabla_z^2 \hat{f}_n(z,x_0) = \sum_{i=1}^n w_n(x_i,x_0) \nabla_z^2 F(z;y_i)$,
respectively, where $\nabla_z F(z;y_i) = 4 A (z-By_i)^{3}$ and $\nabla_z^2 F(z;y_i) = 12A (z-By_i)^{2}$ with $A = \diag(a^1,\ldots,a^{d_y})$ and $B = \diag(b^1,\ldots,b^{d_y})$.
We check the Hessian $\nabla_z^2 \hat{f}_n(z,x_0)$, making sure that its smallest eigenvalue is bounded from below and above by some positive numbers in a neighborhood of $\hat{z}_n(x_0)$.
With careful selection of the initial solution $z^{(0)}_n(x_0)$,
the algorithm can achieve quadratic convergence.

\begin{figure}[ht]
\FIGURE{\includegraphics[width=\textwidth]{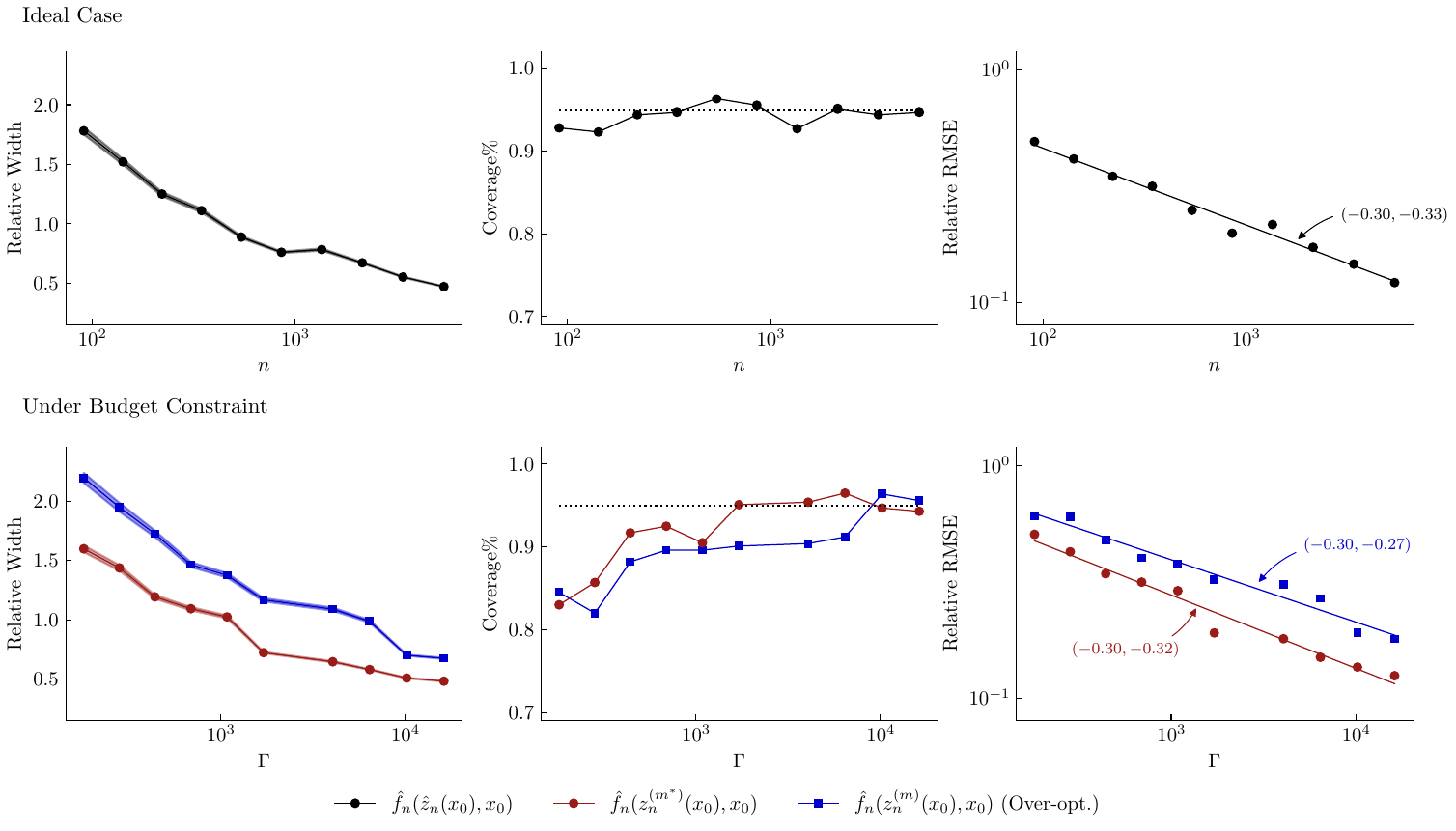}}
{High-order Polynomial Function. \label{fig:exp-Syn-0.25}}
{The first number in each pair of parentheses shows the theoretical convergence rate of relative RMSE: $n^{-(1-\delta d_x)/2}$ for $\hat{f}_n(\hat{z}_n(x_0),x_0)$ (Theorem~\ref{thm:clt-cont}), and $\Gamma^{-(1-\delta d_x)/2}$ for both $\hat{f}_n(z^{(m^*)}_{n}(x_0),x_0)$ (Theorem~\ref{thm:superlin}) and $\hat{f}_n(z^{(m)}_{n}(x_0),x_0)$ (Theorem~\ref{thm:superlin-overopt}), up to logarithmic factors. The second number indicates the empirical slope obtained from regressing log relative RMSEs on $\log n$ or $\log \Gamma$.  ``N/A'' denotes cases where the theoretical rate is not applicable. The optimization algorithm is superlinearly convergent.}
\end{figure}

In Figure~\ref{fig:exp-Syn-0.25}, we evaluate the performance of the confidence intervals~\eqref{eq:ci-wSAA} and \eqref{eq:ci-wSAA-constrained}. For the budget-constrained wSAA estimator, we compare two budget allocation rules: the optimal allocation with $m^* = \kappa^* \log\log\Gamma$ from Theorem~\ref{thm:superlin} and the over-optimizing strategy with $m = \kappa^* \log\Gamma$ from Theorem~\ref{thm:superlin-overopt}, where $\kappa^* = 1/\log2$.
The results validate our theoretical findings for superlinearly convergent algorithms, particularly the confidence intervals' validity and the wSAA estimator's convergence rate.

\section{Experimental Design Details}\label{app:exp-design}

\subsection{Parameters of Optimization Algorithms}\label{app:alg-para}

Optimization algorithms for solving the wSAA problem~\eqref{eq:cso-wSAA}
follow a general form $z^{(t)} = z^{(t-1)} + \mu_t (\Pi^{(t-1)}_{\scrZ} (z^{(t-1)} + g(z^{(t-1)})) - z^{(t-1)}), \forall t \in [m]$,
where $\mu_t$ is the stepsize, $g$ is a function, and $\Pi^{(t)}_{\scrZ}$ is an operator that projects iterates onto the feasible region $\scrZ$,
all of which vary by algorithm. The algorithms implemented in our experiments are summarized as follows.

\begin{enumerate}[label=(\roman*)]
    \item Subgradient descent: $\mu_t \equiv \mu_0/\sqrt{m+1}$ for some $\mu_0>0$; $g(z) = -\partial_z \hat{f}_n(z,x_0)$; and $\Pi^{(t)}_{\scrZ}(z) = \argmin_{z' \in \scrZ} \|z'-z\|^2/2$.
    \item Gradient descent: $\mu_t$ is determined by backtracking line search, shrinking by a factor $b$ until satisfying the Armijo condition $\mu_t = \argmax_{\ell \in \N_+} \{ b^\ell: \hat{f}_n(z^{(t)},x_0) \le \hat{f}_n(z^{(t-1)},x_0) + a b^\ell \langle \nabla_z \hat{f}_n(z^{(t-1)},x_0), \Pi_{\scrZ}^{(t)}(z^{(t-1)} + g(z^{(t-1)})) - z^{(t-1)} \rangle \}, \forall t \in [m]$, where $a \in (0,0.5)$ and $b \in (0,1)$ are constants; $g(z) = -\nabla_z \hat{f}_n(z,x_0)$; and $\Pi^{(t)}_{\scrZ}(z) = \argmin_{z' \in \scrZ} \|z'-z\|^2/2$.
    \item Newton's method: $\mu_t$ is determined by the same line search method; $g(z) = -\nabla_z^2 \hat{f}_n(z,x_0)^{-1} \nabla_z \hat{f}_n(z,x_0)$; and $\Pi_{\scrZ}^{(t)}(z) = \argmin_{z' \in \scrZ} \|z'-z\|^2_{H^{(t)}}/2$ with $H^{(t)} \coloneqq \nabla_z^2 \hat{f}_n(z^{(t)},x_0)$, where $\|z\|_{H} \coloneqq \sqrt{z^\intercal H z}$. (This projection operator ensures that the Armijo rule is enforced along the feasible directions.)
\end{enumerate}

\subsection{Cross-validation}\label{app:CV}

We use $k$-fold cross-validation to select tuning parameters.
These may include the bandwidth constant $h_0$, the stepsize constant $\mu_0$, and the initial solution $z^{(0)}_n(x_0)$, depending on the specific experiment.
Let $\Xi$ denote this set of tuning parameters.
We divide the dataset $\scrD_n$ into $k$ folds.
Let $\{\scrI_\ell: \ell \in [k]\}$ denote a collection of $k$ equal-sized partitions of $[n]$ (assuming $n$ is divisible by $k$ for simplicity). Let $\scrD_{-\ell}$ denote the set of data from $\scrD_n$ with the $\ell$-th fold $\{(x_i,y_i)\}_{i \in \scrI_\ell}$ excluded.
The wSAA problem on the dataset $\scrD_{-\ell}$ is formulated as $\min_{z \in \scrZ} \left\{ \hat{f}_{-\ell}(z, x_0) \coloneqq \sum_{i \in [n] \setminus \scrI_\ell} w_{-\ell}(x_i,x_0) F(z;y_i) \right\}$,
where $w_{-\ell}(x_i,x_0) \coloneqq \frac{K((x_i-x_0)/h_n)}{\sum_{i \in [n] \backslash \scrI_\ell} K((x_i-x_0)/h_n)}$ for all $i \in [n] \backslash \scrI_\ell$.
Let $\hat{z}_{-\ell}(x_0;\Xi)$ denote the optimal solution to this problem.
For experiments without computational budget constraints, we only tune $\Xi = \{h_0\}$. This parameter is selected based on the average out-of-sample performance across $k$ folds, given by $\mathsf{CV}_k(\Xi) = \frac{1}{k} \sum_{\ell \in [k]} \sum_{i \in \scrI_\ell}F(\hat{z}_{-\ell}(x_0;\Xi);y_i)$.
Moreover, for experiments with computational budget constraints,
we let $z^{(m)}_{-\ell}(x_0;\Xi)$ denote the solution after $m$ iterations from the initial point  $z^{(0)}_n(x_0)$.
The tuning parameters are $\Xi = \{h_0,\mu_0,z^{(0)}_n(x_0)\}$ for subgradient descent, and $\Xi = \{h_0,z^{(0)}_n(x_0)\}$ for gradient descent and Newton's method.
The average out-of-sample performance across $k$ folds then becomes $\mathsf{CV}_k(\Xi) = \frac{1}{k} \sum_{\ell \in [k]} \sum_{i \in \scrI_\ell}F(z^{(m)}_{-\ell}(x_0;\Xi);y_i)$.

\subsection{Training of CWGAN Simulator} \label{app:CWGAN}

The loss function of the CWGAN with gradient penalty is defined as $\scrL_{\text{CWGAN-gp}}(\theta_{g}, \theta_{c}) =\E_{Y|X} \left[C(Y|X;\theta_{c}) \right] - \E_{U|X} \left[ C(G(U|X;\theta_{g})|X;\theta_{c}) \right] - \lambda_{\text{gp}} \E_{Y,U|X}[(\|\nabla_{\tilde{Y}} C(\tilde{Y}|X;\theta_{c})\|-1)^+]$.
Here, $G(\cdot;\theta_{g})$ is the generator producing samples with latent noise $U$ (e.g., drawn from standard normal or uniform distribution) to mimic the true conditional distribution of $Y$ given $X$. $C(\cdot;\theta_{c})$ is the critic distinguishing between the generated sample $\hat{Y} = G(U|X;\theta_{g})$ and the real sample $Y$ for the given $X$.
The interpolated sample $\tilde{Y} = \epsilon Y + (1-\epsilon) \hat{Y}$ with $\epsilon \sim \textsf{Unif}[0,1]$ lies on the segment joining $Y$ and $\hat{Y}$.
The gradient penalty term enforces a soft 1-Lipschitz constraint on the critic by penalizing deviations of its gradient norms from one.
Inspired by the conditional GAN literature, we augment $\scrL_{\text{CWGAN-gp}}(\theta_{g}, \theta_{c})$ with an additional reconstruction loss term to enhance training stability. The reconstruction loss is typically in the form of $\ell_p$ norm; in our implementation, we set $p = 2$. It is defined as $\scrL_{\text{rec}}(\theta_g) = \lambda_{\text{rec}} \E_{Y,U|X}[\|Y-G(U|X;\theta_{g})\|]$, where $\lambda_{\mathrm{rec}}>0$ is a hyperparameter that balances between the adversarial objective $\scrL_{\text{CWGAN-gp}}$ and this reconstruction objective.
The generator is therefore trained to both deceive the critic and reproduce the ground truth.
We optimize $\theta_g$ and $\theta_c$---the parameters of generator and critic---using the Optimistic Adam optimizer \citep{daskalakis2018training}. It alleviates the phenomenon of limit cycling and improves convergence in the adversarial, min–max game; that is, {$\min_{\theta_{g}} \max_{\theta_{c}} \ \scrL_{\text{CWGAN-gp}} + \lambda_{\text{rec}} \scrL_{rec}(\theta_g)$.
In our implementation, CWGAN is trained with a batch size of 512 (approximately 5\% of the training data) for up to 10,000 epochs. The generated hourly demand values are restricted to the range $[0,10000]$.
To achieve optimal model performance, we apply the Asynchronous Successive Halving Algorithm \citep{Li20}---implemented via the framework Ray (\url{https://www.ray.io/})---to tune the architecture of neural nets as well as other key hyperparameters such as the learning rate.

\end{document}